\documentclass{amsart}
\usepackage{amssymb}
\usepackage{xcolor}
\usepackage{extarrows}
\usepackage{thmtools,mathtools}
\usepackage{hyperref}
\usepackage[nameinlink]{cleveref}
\usepackage{mleftright}
\definecolor{royalblue}{rgb}{0.0, 0.14, 0.9}
\definecolor{applegreen}{rgb}{0.55, 0.71, 0.0}
\definecolor{webbrown}{rgb}{0.65, 0.56, 0.16}
\hypersetup{%
    colorlinks=true, linktocpage=true, pdfstartview={FitV null},%
    breaklinks=true, pageanchor=true, pdfpagemode=UseOutlines,
    plainpages=false, bookmarksnumbered, bookmarksopen=true, bookmarksopenlevel=1,
    hypertexnames=false, pdfhighlight=/O,
    urlcolor=webbrown, linkcolor=royalblue, citecolor=applegreen,
    pdftitle={An equational axiomatisation of metrically complete lattice-ordered groups},%
    pdfauthor={Marco Abbadini, Vincenzo Marra, and Luca Spada},%
    pdfsubject={Equational axiomatisation of a variety equivalent to metrically complete Archimedean unital lattice-ordered groups and metrically complete MV-algebras},%
    pdfkeywords={lattice-ordered groups, unital lattice-ordered groups, Archimedean lattice-ordered groups, metrically complete lattice-ordered groups, MV-algebras, Mundici equivalence, metric completeness, infinitary algebras, infinitary varieties, equational axiomatisation, Linton-Lawvere theory, normal a-spaces},%
    pdfcreator={pdfLaTeX},%
}   
    
\usepackage{tikz-cd}
\usepackage{bm}
\usepackage[italian,backgroundcolor=blue!5!white,textsize=small]{todonotes}
\usepackage{enumitem}
\usepackage{xargs}
\usepackage{tikz}
\usepackage{pgfplots}
 \pgfplotsset{compat=1.16}
\usetikzlibrary{shapes.misc}
\usepackage{subcaption}
\declaretheoremstyle[
spaceabove=6pt, 
spacebelow=6pt,
headfont=\normalfont\bfseries,
notefont=\mdseries, 
notebraces={(}{)},
bodyfont=\itshape,
postheadspace=.5em,
]{theorem}
\declaretheorem[numberwithin=section,style=theorem,refname={Theorem,Theorems},  Refname={Theorem,Theorems}, name=Theorem]{theorem}
\declaretheorem[sibling=theorem,style=remark,refname={Remark,Remarks}, Refname={Remark,Remarks},name=Remark]{remark}
\declaretheorem[sibling=theorem,style=theorem,refname={Proposition,Propositions},  Refname={Proposition,Propositions},name=Proposition]{proposition}
\declaretheorem[sibling=theorem,style=theorem, name=Lemma, refname={Lemma,Lemmas},  Refname={Lemma,Lemmas}]{lemma}
\declaretheorem[sibling=theorem,style=theorem,refname={Corollary,Corollaries},  Refname={Corollary,Corollaries},name=Corollary]{corollary}

\declaretheorem[sibling=theorem,style=definition,refname={Notation,Notations},  Refname={Notation,Notations},name=Notation]{notation}
\declaretheorem[sibling=theorem,style=definition,refname={Definition,Definitions},  Refname={Definition,Definitions},name=Definition]{definition}
\declaretheorem[sibling=theorem,style=definition,refname={Example,Examples},  Refname={Example,Examples},name=Example]{example}

\numberwithin{equation}{section}

\renewcommand{\leq}{\leqslant}
\renewcommand{\le}{\leqslant}
\renewcommand{\geq}{\geqslant}
\renewcommand{\ge}{\geqslant}
\renewcommand{\preceq}{\preccurlyeq}

\newcommand{\df}{\coloneqq}
\newcommand{\adj}{\raisebox{.07\baselineskip}{\ensuremath{\scriptscriptstyle \top}}}
\newcommand{\KHa}{\mathsf{KH}_{\mathsf{a}}}
\newcommand{\ISP}{\mathbb{ISP}}

\newcommand{\succaux}[3]{\ensuremath{\langle#1\mid{#2\in#3}\rangle}}
\newcommandx*{\succustom}[3][2=n, 3 = \Npos]{\succaux{#1}{#2}{#3}}

\renewcommandx*{\succ}[3][2=n, 3 = \Npos]{\bm{#1}}
\newcommand{\asucc}[1]{\bm{#1}} 
\newcommandx*{\init}[2][2=n]{\ensuremath{\bm{#1}_{\leq#2}}}
\newcommand{\shift}[1]{\operatorname{\uparrow}{#1}}

\newcommand{\cat}[1]{\mathsf{#1}}
\newcommand{\Set}{\mathsf{Set}}
\newcommand{\G}{\mathsf{G}} 
\newcommand{\cG}{\overline{\G}} 
\newcommand{\KH}{\mathsf{KH}} 
\newcommand{\CMV}{\mathsf{CMV}}
\newcommand{\MV}{\mathsf{MV}}
\newcommand{\cMV}{\overline{\MV}} 
\newcommand{\op}{\mathrm{op}}

\DeclareMathOperator{\U}{U}
\DeclareMathOperator{\F}{F}

\newcommand{\R}{\mathbb{R}}
\newcommand{\Q}{\mathbb{Q}}
\newcommand{\Z}{\mathbb{Z}}
\newcommand{\Nzero}{\mathbb{N}_0}
\newcommand{\Npos}{\mathbb{N}}

\DeclareMathOperator{\den}{den}
\DeclareMathOperator{\lcm}{lcm}
\DeclareMathOperator{\dl}{\mathtt{d}_{\ell}}
\DeclareMathOperator{\dMV}{\mathtt{d}_{\mathrm{MV}}}
\newcommand{\norml}[1]{\lVert#1\rVert_{\ell}}
\newcommand{\normMV}[1]{\lVert#1\rVert_{\mathrm{MV}}}
\newcommand{\mvm}[1]{\left[#1\right]}
\newcommand{\lm}[1]{#1}
\DeclareMathOperator{\cha}{c} 
\newcommand{\limn}{\underset{n\to+\infty}{\lim}}

\title[An equational axiomatisation of metrically complete $\ell$-groups]{An equational axiomatisation of metrically complete lattice-ordered groups}

\author{Marco Abbadini}
\address[Marco Abbadini]{Research Institute in Mathematics and Physics (IRMP), Université catholique de Louvain, 
Chem.\ du Cyclotron 2, 1348 Louvain-la-Neuve, Belgium}

\author{Vincenzo Marra}
\address[Vincenzo Marra]{Dipartimento di Matematica `Federigo Enriques', Universit\`a degli Studi di Milano, via Cesare Saldini 50, 20133 Milano, Italy}

\author{Luca Spada}
\address[Luca Spada]{Dipartimento di Matematica, Universit\`a degli Studi di Salerno, Piazza Renato Caccioppoli, 2, 84084 Fisciano (SA), Italy}

\begin{document}

\begin{abstract}
    We provide an equational axiomatisation of a variety equivalent to the category of Archimedean unital lattice-ordered groups that are metrically complete with respect to the metric induced by the order unit. While we use an operation of countably infinite arity, we do provide a finite equational axiomatisation. The countable bound on the arity is optimal: we show that the resulting category is not equivalent to any finitary variety---nor, indeed, to any first-order definable class of structures in a finitary language, as well as to any class of finitary algebras closed under isomorphisms, subalgebras and products.
\end{abstract}

\keywords{Lattice-ordered group, unital $\ell$-group, MV-algebra, Mundici equivalence, metric completeness, infinitary variety, equational axiomatisation, Linton--Lawvere theory, normal a-spaces}

\subjclass[2020]{Primary: 06F20. Secondary: 06D35, 03G20, 08B05, 18C05, 08A65}

\maketitle

\tableofcontents

\section*{Introduction} \label{s:intro}
Stone duality asserts that the category of compact Hausdorff zero-dimensional spaces is dually equivalent to the finitary variety of Boolean algebras.  Beyond dimension zero, the situation changes substantially.  It has been known since the work of Duskin \cite{Duskin} and Negrepontis \cite{Negrepontis1969,Negrepontis1971} that the opposite of the category of compact Hausdorff spaces is monadic over sets, and hence equivalent to a possibly infinitary variety.  Isbell \cite{Isbell} subsequently showed that an algebraic theory may be generated by finitely many finitary operations together with an operation of countably infinite arity.  A finite equational axiomatisation of this theory was later obtained by Marra and Reggio in \cite{MarraReggio}. The countable bound on the arity is provably optimal, in that the category of compact Hausdorff spaces is not dually equivalent to any finitary variety \cite[Prop.~1.3]{MarraReggio}.

In light of Yosida's version of Stone--Gelfand duality \cite{Yosida1941}, the axiomatisation in \cite{MarraReggio} is, up to categorical equivalence, an equational description of metrically complete unital vector lattices.  From this point of view, the primitive operation of countably infinite arity captures metric-completeness: indeed, on algebras of continuous functions it computes the limit of a suitably controlled sequence. We refer the reader to the introductions of \cite{AbbadiniMarraEtAl2025, MarraReggio} for fuller historical background, and for a discussion of the significance of these results in algebraic logic, topology, and functional analysis.

The present paper weakens the real vector space structure from Yosida's setting to the Abelian group structure, thereby addressing the analogous axiomatisation problem for Archimedean unital lattice-ordered groups that are complete in the metric induced by their (order) unit.  We call these objects \emph{metrically complete $\ell$-groups}. Our main result is that the category $\cG$ of metrically complete $\ell$-groups and unital $\ell$-group homomorphisms is equivalent to a variety of infinitary algebras admitting a finite equational axiomatisation.  The language consists of the usual finitary operations of MV-algebras \cite{MundiciAlg} together with a single operation of countably infinite arity.

The r\^{o}le of MV-algebras is of the essence here.  Under Mundici's equivalence \cite{M86} between unital Abelian lattice-ordered groups and MV-algebras, Archimedean unital lattice-ordered groups correspond precisely to semisimple MV-algebras.  Moreover, the metric induced by the unit restricts to the natural metric on the corresponding MV-algebra, and completeness is preserved and reflected by the equivalence.  
We call an MV-algebra \emph{metrically complete} when it is semisimple and complete in its metric.
Consequently, the category $\cG$ of metrically complete $\ell$-groups and unital $\ell$-group homomorphisms is equivalent to the category $\cMV$ of metrically complete MV-algebras and MV-homomorphisms.

It is on this basis that we then construct an infinitary variety $\CMV$.  Its objects are MV-algebras equipped with one operation
\[
\gamma\colon A^{\Npos}\longrightarrow A
\]
of countably infinite arity.  Informally, $\gamma$ evaluates the limit of a sequence.  However, since an arbitrary sequence need not be Cauchy, the operation cannot simply be interpreted as its limit.
We therefore define, for each positive integer $n$, an MV-term
\[\rho_n(x_1,\dots,x_n)\]
with the following two properties.
\begin{enumerate}
    \item For every sequence $(x_n)_{n\in \Npos}$ in a semisimple MV-algebra, the sequence
\[(x_n')_{n \in \Npos} \coloneqq\bigl(\rho_n(x_1,\dots,x_n)\bigr)_{n\in \Npos} \]
is \emph{fast}, i.e., the distance between $x_n'$ and $x_{n+1}'$ is less than or equal to $\frac{1}{2^{n+1}}$.
    \item If the original sequence $(x_n)_{n\in \Npos}$ is already fast, then, for every $n$,
    \[\rho_n(x_1,\dots,x_n)=x_n.\] 
\end{enumerate}
Since a metric space is complete precisely when all fast sequences converge, in a metrically complete MV-algebra we may set
\[\gamma\bigl((x_n)_{n\in \Npos}\bigr)=\lim_{n\to\infty}\rho_n(x_1,\dots,x_n). \]

The axioms of $\CMV$ ensure that this is exactly the meaning of $\gamma$.  Two of these axioms force the underlying MV-algebra to be semisimple, while the remaining ones prescribe that $\gamma\bigl((x_n)_{n\in \Npos}\bigr)$ approximates $\rho_n(x_1,\dots,x_n)$ with error at most $\frac{1}{2^n}$; it follows that every $\CMV$-algebra is metrically complete.
Conversely, we show that every metrically complete MV-algebra admits a unique expansion to a $\CMV$-algebra. 
Furthermore, every MV-homomorphism between metrically complete MV-algebras preserves the terms $\rho_n$ and, being non-expansive, also preserves the limits of the resulting sequences; consequently, it preserves $\gamma$.
Therefore, the forgetful functor
\[\U\colon\CMV\longrightarrow\MV\]
is full and faithful, and its image consists precisely of the metrically complete MV-algebras.  We thus obtain equivalences
\[\CMV \cong \cMV \simeq \cG.\]

The  axiomatisation we obtain in this manner contains an infinite set of equations, one for each prescribed approximation bound. With a further argument, we reduce it to an equivalent but \emph{finite} set of equations. 
Thus, the infinitary character of $\CMV$ does not lie in the number of equational axioms, but rather in the arity of $\gamma$.

Having obtained a variety equivalent to $\cMV$, one may speak of free metrically complete MV-algebras: these are the objects corresponding, through the above equivalence, to free $\CMV$-algebras.  We prove a representation theorem for these free objects, as follows. We first prove that the standard MV-algebra $[0,1]$, equipped with the operation $\gamma$, generates $\CMV$ not only as a variety but also under just isomorphisms, subalgebras, and products:
\[\CMV=\ISP([0,1]).\]
The operations of $\CMV$ can consequently be identified with suitable functions on powers of $[0,1]$.  Using the Stone--Weierstrass theorem for arithmetic spaces developed in \cite{AbbadiniMarraEtAl2025}, we prove that the interpretations in $[0,1]$ of the $\CMV$-terms of arity $\lambda$ are exactly the \emph{a-maps} from $[0,1]^\lambda$ to $[0,1]$; these are, by definition, the continuous functions 
\[[0,1]^\lambda \to [0,1]\]
that \emph{decrease} (with respect to the divisibility order) \emph{denominators} (for the appropriate notion of denominators for elements of Tychonoff cubes). It follows that the free $\CMV$-algebra on a set $I$ is the algebra of a-maps $[0,1]^I \to [0,1]$. In terms of algebraic theories, these results say that $\CMV$ is the category of algebras of the infinitary Linton--Lawvere theory whose objects are powers of $[0,1]$ and whose morphisms are a-maps between them \cite{Lawvere,Linton}.

The use of an operation of infinite arity in our axiomatisation is unavoidable.  We prove, in fact, a stronger negative result: the categories $\cG$ and $\cMV$ of metrically complete $\ell$-groups and metrically complete MV-algebras are not equivalent to any elementary class of structures in a finitary language, with all homomorphisms as morphisms.  In particular, they are not equivalent to any finitary variety.  The proof uses the duality between metrically complete $\ell$-groups and normal a-spaces obtained in \cite{AbbadiniMarraEtAl2025} and an obstruction to finite concreteness originating in the opposite of the category $\KH$ of compact Hausdorff spaces: no faithful functor
\[\KH^{\op}\longrightarrow\Set\]
preserves directed colimits \cite{LiebermanRosickyEtAl2023,Abbadini2026}.  We transfer this obstruction to the opposite of the category of normal a-spaces, and hence to $\CMV$.  Since the ordinary underlying-set functor of an elementary class of finitary structures preserves directed colimits \cite{Richter1971}, the claimed impossibility follows.  Thus, the countable arity of the operations in our axiomatisation is optimal.

As a further negative result, we prove that $\cG$ and $\cMV$ are not equivalent to any $\ISP$-class of finitary algebras, i.e., to any class of finitary algebras closed under isomorphisms, subalgebras and products.
While this statement a priori requires consideration only of finitary algebras, our main result presenting $\cMV$ as a variety of infinitary algebras plays a key role, as it ensures a purely categorical property---namely Barr-exactness---that translates to closure under homomorphic images in any $\ISP$-class of algebras.

\section{Preliminaries on MV-algebras}

In this paper, $\Npos\df\{1,2,3,\dots\}$ and $\Nzero\df\Npos\cup\{0\}$.

Two main ingredients towards an equational axiomatisation of the class of metrically complete $\ell$-groups are MV-algebras and Mundici's functor $\Gamma$, which we briefly recall here; standard references are \cite{MundiciAlg,MundiciAdv}.
\begin{definition}[MV-algebras] \label{d:MV-algebra}
An \emph{MV-algebra} is a structure $(A,\oplus,\neg,0)$ with the following properties.
\begin{align}
\label{d:MV-algebras:item1}\tag{{MV} 1} &(A,\oplus,0) \text{ is a commutative monoid}.\\
\label{d:MV-algebras:item2}\tag{{MV} 2} &x \oplus \neg0 = \neg0.\\
\label{d:MV-algebras:item3}\tag{{MV} 3} &\neg(\neg x) = x.\\
\label{d:MV-algebras:item4}\tag{{MV} 4} &\neg(\neg x \oplus y) \oplus y = \neg(\neg y\oplus x) \oplus x.
\end{align}
\end{definition}

It is customary to define the following derived operations in the language of MV-algebras:
\begin{align*}
&1\df \neg 0,\quad x \odot y\df \neg(\neg x \oplus \neg y), \quad x \ominus y\df  x \odot \neg y, \quad\text{and,} \\
&\text{for any }n \in \Nzero,\quad
[n]x\coloneqq \underbrace{x\oplus\cdots\oplus x}_{n\ \mathrm{times}},
\end{align*}
where, by convention, $[0]x = 0$.
Every MV-algebra can be endowed with the structure of a distributive lattice by defining 
\[
x \vee y\df(x \ominus y) \oplus y\quad\text{and}\quad x \wedge y\df \neg (\neg x \vee \neg y).
\]
The lattice order is obtained as usual as
\[
x \leq y\iff x \wedge y = x \iff x \vee y = y.
\]
Hence, inequalities can be defined by equations; we will use this fact without explicitly mentioning it.
Moreover, we have the following derived operation:
\[
\cha(x,y) \df (x \ominus y) \vee (y \ominus x).
\]
The element $\cha(x,y)$ is called the \emph{Chang distance} between $x$ and $y$. 

\begin{example}[{The standard MV-algebra $[0,1]$}] \label{e:standard-MV-algebra}
The \emph{standard MV-algebra} is the real unit interval $[0,1]$ with operations $x \oplus y = \min(x+y,1)$ and $\neg x = 1-x$.  Easy calculations show that, in this algebra,
\begin{align*}
&x \odot y = \max \{0,x+y-1\}, \quad x \ominus y = \max \{0,x-y\},\quad \cha(x,y) = \lvert x-y \rvert,\\ 
&x\vee y = \max \{x,y\}, \quad x\wedge y = \min \{x,y\}.
\end{align*}
We will always refer to $[0,1]$ endowed with this canonical structure of an MV-algebra.
\end{example}

Recall that an $\ell$-group is an algebraic structure $(G,+,-,0,\wedge,\vee)$ such that:
\begin{enumerate}
    \item\label{l-group:axiom1} \((G,+,-,0)\) is a group;
    \item\label{l-group:axiom2} \((G,\vee,\wedge)\) is a lattice;
    \item\label{l-group:axiom3} for all \(x,y,z \in G\),
    \[
z+(x\vee y)=(z+x)\vee (z+y),\qquad z+(x\wedge y)=(z+x)\wedge (z+y),
\]
and similarly on the right.
\end{enumerate}
The condition \eqref{l-group:axiom3} can be equivalently expressed by saying that the group operation is order-preserving: for all \(x,y,z\in G\),
    \[
    x \leq y \text{ implies } z+x \leq z+y \quad\text{and}\quad x+z \leq y+z,
    \]
    where the order \(\leq\) is the lattice order defined by
    \[
    x \leq y \iff x \wedge y = x \iff x \lor y = y.
    \]
An $\ell$-group is called \emph{Abelian} if its underlying group is Abelian.

\begin{definition}[Strong order unit]\label{d:strong-order-unit}
Let \(G\) be an \(\ell\)-group. An element \(1\in G\) with \(1\geq 0\) is called a \emph{(strong) order unit} if for every \(x\in G\) there is \(n\in\mathbb{N}\)  such that
\[
-n1 \leq x \leq n1.
\]
\end{definition}
A \emph{unital} $\ell$-group is a pair $(G,1)$ consisting of an $\ell$-group $G$ and a designated order unit $1\in G$. In the category of unital $\ell$-groups, the morphisms are the $\ell$-group homomorphisms preserving the designated order unit.

\begin{definition}[Mundici's functor $\Gamma$]\label{d:gamma}
The functor $\Gamma$ from the category of unital Abelian $\ell$-groups to the category of MV-algebras is defined as follows.
If $(G,1)$ is a unital Abelian $\ell$-group, then $\Gamma(G,1) \df \{g \in G\mid 0\leq g \leq 1\}$, where the MV-algebraic operations are defined as follows:
\[
x\oplus y\df (x+y) \wedge 1\qquad\text{ and }\qquad \neg x\df 1-x.
\]
If $h$ is a unit-preserving homomorphism between two unital Abelian $\ell$-groups $(G,1)$ and $(H,1)$, then $\Gamma(h)$ is defined as the restriction of $h$ to $\Gamma(G,1)$ and $\Gamma(H,1)$.  The fact that the homomorphism $h$ preserves the unit guarantees that $\Gamma(h)$ ranges in $\Gamma(H,1)$.
\end{definition}
\begin{theorem}[{Mundici's equivalence \cite[Thm.~3.9]{M86}}] \label{th:gamma}
The above-defined functor $\Gamma$ has a quasi-inverse $\Xi$, i.e., $\Gamma$ and  $\Xi$ form a categorical equivalence.
\end{theorem}

As an instance of a general definition in algebra, an MV-algebra is called \emph{semisimple} if it is a subdirect product of simple MV-algebras. Semisimple MV-algebras can be described as algebras of $[0,1]$-valued functions:
\begin{proposition}[{\cite[Prop.~3.6.1]{MundiciAlg}}]\label{p:semisimple-algebras-of-functions}
An MV-algebra is semisimple if and only if it is isomorphic to a subalgebra of $[0,1]^X$ for some set $X$, where MV-operations are calculated pointwise.
\end{proposition}

\begin{lemma}\label{l:MV-subd-lin-ord}
    Every MV-algebra is a subdirect product of linearly ordered MV-algebras. 
\end{lemma}
\begin{proof}
    This was originally proved in \cite[Lem.~3]{Chang-completeness}.  See \cite[Thm.~1.3.3]{MundiciAlg} for a textbook reference.
\end{proof}

\begin{theorem}[Di Nola's Representation Theorem]\label{t:dinola}
    Every non-trivial linearly ordered MV-algebra embeds into some ultrapower of the standard MV-algebra $[0,1]$. As a consequence, $[0,1]$ generates the class of MV-algebras as a quasi-variety.
\end{theorem}
\begin{proof}
    For a textbook reference of the fact that every non-trivial linearly ordered MV-algebra embeds into some ultrapower of the standard MV-algebra $[0,1]$, see the proof of \cite[Thm.~9.5.1]{MundiciAlg}.     This was originally proved in \cite{di1991representation} and strengthened in \cite{di2010representation}.

    The fact that $[0,1]$ generates the class of MV-algebras as a quasi-variety follows from this fact together with \cref{l:MV-subd-lin-ord}.
\end{proof}

\begin{lemma} \label{l:MV-basic-properties}
Let $A$ be an MV-algebra.
\begin{enumerate}
\item\label{l:MV-basic-properties:item-monotone} The operation $\oplus$ is monotone on both coordinates, and the operation $\ominus$ is monotone on the first coordinate and antitone on the second one.
\item\label{l:MV-basic-properties:item-neg-antitone} The operation $\neg$ is antitone.
\item\label{l:MV-basic-properties:item-cancellation-identity} For all $x,y \in A$, $(x\oplus (y\ominus x)) \ominus (x\ominus y) = y$.
\end{enumerate}
\end{lemma}

\begin{proof}
Items \eqref{l:MV-basic-properties:item-monotone} and
\eqref{l:MV-basic-properties:item-neg-antitone} are standard; see, e.g.,
\cite{M86}.

For \eqref{l:MV-basic-properties:item-cancellation-identity}, since the identity
is an equation in the language of MV-algebras, by \Cref{t:dinola} it is enough
to check it in the standard MV-algebra $[0,1]$. Let $x,y\in[0,1]$. If $y\geq x$,
then $y\ominus x=y-x$ and $x\ominus y=0$, and so
\[
    (x\oplus (y\ominus x))\ominus (x\ominus y)
    =(x\oplus (y-x))\ominus 0
    =y\ominus 0
    =y.
\]
If $x>y$, then $y\ominus x=0$ and $x\ominus y=x-y$, and so
\[
    (x\oplus (y\ominus x))\ominus (x\ominus y)
    =(x\oplus 0)\ominus (x-y)
    =x\ominus (x-y)
    =y.
\]
Thus the identity holds in $[0,1]$, and hence in every MV-algebra.
\end{proof}

Let \(G\) be an \(\ell\)-group and \(x,y \in G\) with $y \geq 0$. We say that $x$ is \emph{infinitesimal with respect to} $y$ if
\[
-y \leq \lm{n}x \leq y \text{ for all } n\in\mathbb{N}.
\]
The $\ell$-group \(G\) is said to be \emph{Archimedean} if, for every $y\in G$ with $y \geq 0$, the only infinitesimal element of $G$ with respect to $y$ is $0$.
Every Archimedean $\ell$-group is Abelian \cite[Thm.~11.1.3]{BigardKeimelEtAl1977}. Therefore, in the following, we omit this qualification when the group is assumed to be Archimedean.
\begin{remark}\label{r:unital-simplifies-Archimedean}
Using the definition of strong order unit, it is easy to see that a unital Abelian $\ell$-group is Archimedean if and only if there are no non-zero infinitesimal elements with respect to the unit.
\end{remark}

The class of MV-algebras corresponding to Archimedean unital $\ell$-groups can be characterised in universal algebraic terms.  
\begin{lemma} \label{l:gamma:Archimedean}
The functor $\Gamma$ of \Cref{th:gamma} restricts to an equivalence between Archimedean unital $\ell$-groups and semisimple MV-algebras.
\end{lemma}
\begin{proof}
It is enough to check that the functor $\Gamma$ sends Archimedean unital $\ell$-groups into semisimple MV-algebras and vice versa for $\Xi$.

We begin by checking that $\Gamma$ sends Archimedean unital $\ell$-groups to semisimple MV-algebras.
Let $(G,1)$ be an Archimedean unital $\ell$-group. By the representation theorem for Archimedean $\ell$-groups with (strong) unit (see, e.g., \cite[Cor.~13.2.6]{BigardKeimelEtAl1977}), there are a set $X$ and an $\ell$-embedding $\iota\colon G\hookrightarrow \R^{X}$ sending $1$ to the constant function $1$.
Applying $\Gamma$ and using the definition of $\Gamma(G,1)$ as the unit interval $\{g\in G\mid 0\le g\le 1\}$, we obtain an MV-embedding
\[
\Gamma(\iota)\colon \Gamma(G,1)\hookrightarrow [0,1]^X,
\]
where operations are computed pointwise. By \Cref{p:semisimple-algebras-of-functions} $\Gamma(G,1)$ is semisimple.

Next, we check that $\Xi$ sends semisimple MV-algebras to Archimedean unital $\ell$-groups.
Let $A$ be a semisimple MV-algebra. By \cref{p:semisimple-algebras-of-functions}, there are a set $X$ and an MV-embedding $e\colon A\hookrightarrow [0,1]^X$.
It is well-known that the functor $\Xi$ preserves the injectivity of the underlying functions of morphisms\footnote{See, for example, \cite[Prop.~3.7]{AbbadiniJipsenEtAl2022}, which is in a more general setting.}. Thus, applying $\Xi$, we obtain a unital $\ell$-embedding
\[
\Xi(e)\colon \Xi(A)\hookrightarrow \Xi([0,1]^X).
\]
Since $\Xi([0,1]^X)$ is the unital $\ell$-group of bounded functions from $X$ into $\R$ (with pointwise order), $\Xi([0,1]^X)$ is Archimedean.
Thus, any $\ell$-subgroup of $\Xi([0,1]^X)$ is Archimedean, and so $\Xi(A)$ is Archimedean.
\end{proof}

\begin{theorem}\label{t:quasivariety-generated-by-linear}
Each non-trivial linearly ordered Abelian group generates the class of Abelian $\ell$-groups as a quasi-variety.
\end{theorem}
\begin{proof}
Let $A$ be a non-trivial linearly ordered Abelian $\ell$-group, and $q$ a quasi-equation that holds in $A$.
Let $H$ be an Abelian $\ell$-group, and let us prove that $q$ holds in $H$.
By Birkhoff's subdirect representation theorem, $H$ is a subdirect product of subdirectly irreducible Abelian $\ell$-groups.
Every finitely subdirectly irreducible Abelian $\ell$-group (and, in particular, every subdirectly irreducible Abelian $\ell$-group) is totally ordered (\cite[Thm.~23]{Birkhoff1942}) and non-trivial.
Moreover, by a result of Gurevich and Kokorin \cite[Thm.~1]{GurevichKokorin1963}, all non-trivial totally ordered Abelian $\ell$-groups are universally equivalent; in particular, they satisfy the same quasi-equations.
Therefore, $q$ holds in every subdirectly irreducible Abelian $\ell$-group, and hence in $H$.
\end{proof}

\section{Metrically complete MV-algebras}

\begin{definition}\label{d:ell-norm}
Let $(G,1)$ be an Abelian unital $\ell$-group.
\begin{enumerate}
    \item For $g\in G$, we set
    \begin{equation} \label{eq:norml}
        \norml{g} \coloneqq \inf\left\{\frac{p}{q} \mid p \in \Nzero,\, q \in \Npos,\, -p 1 \le q g\le p1\right\}.
    \end{equation}
    We call the function $\norml{\,} \colon G \to \R$ the \emph{pseudonorm} (or simply the \emph{norm} if $G$ is Archimedean) on $G$ induced by the unit $1$. See \cite[Sec.~7, p.~120]{Goodbook}.
    \item 
    For $g,h \in G$, we set 
    \[
    \dl(g,h) \coloneqq \norml{g-h};
    \]
    moreover, we call the function $\dl \colon G \times G \to \R$ the \emph{pseudometric} (or simply the \emph{metric} if $G$ is Archimedean) on $G$ induced by the unit $1$ (see \cite[Sec.~15, p.~237]{Goodbook}). The function $\dl$ is indeed a pseudometric on $G$ (see \cite[Prop.~7.12 and page~237]{Goodbook}).
\end{enumerate}
\end{definition}

Recall that, given an Abelian $\ell$-group $G$, for every $g \in G$ one sets
\[
\lvert g \rvert \coloneqq g \lor (-g).
\]
The element $\lvert g \rvert$ is called the \emph{absolute value} of $g$ (\cite[Def.~1.3.8]{BigardKeimelEtAl1977}); by \cite[1.3.10]{BigardKeimelEtAl1977},
\begin{equation} \label{eq:absolute-value-is-positive}
    \lvert g \rvert \geq 0.
\end{equation}
The condition $-p 1 \le q g \le p1$
appearing in \eqref{eq:norml} is equivalent to $\lvert qg \rvert \leq p1$ (this is an immediate standard fact; see \cite[1.3.9]{BigardKeimelEtAl1977}), which is also equivalent to $q\lvert g \rvert \leq p1$ (\cite[1.3.14]{BigardKeimelEtAl1977}).
Since from $\lvert g \rvert \geq 0$ one deduces $\lvert \lvert g \rvert \rvert = \lvert g \rvert$, we conclude
\begin{equation} \label{eq:norm-is-norm-of-absolute-value}
    \norml{g} = \norml{\lvert g \rvert}.
\end{equation}

The function $\dl$ is a metric if and only if $G$ is Archimedean, as proved in the following lemma.
\begin{lemma}\label{lem:kernel-unit-norm}
Let $(G,1)$ be a unital Abelian $\ell$-group and let $\|\cdot\|_\ell$ be the  pseudonorm on $G$ induced by the unit $1$. For every $g\in G$, $\|g\|_\ell=0$ if and only if for every $n\in\Npos$ one has $n\lvert g\rvert\le 1$.
That is, $\|g\|_\ell=0$ if and only if $g$ is an infinitesimal with respect to $1$. Therefore, the pseudometric $d_\ell(g,h)=\|g-h\|_\ell$ is a metric if and only if $(G,1)$ is Archimedean.
\end{lemma}
\begin{proof}
Assume $\|g\|_\ell=0$ and fix $n\in\Npos$. By definition of infimum, there are $p\in\Nzero$ and
$q\in\Npos$ such that $\frac pq<\frac1n$ and $q \lvert g\rvert\le p1$.
Since $\frac pq<\frac1n$ we have $np<q$, and so $np1\le q1$. Multiplying $q\lvert g\rvert\le p1$ by $n$ gives
$nq\lvert g\rvert\le np1\le q1$. As $q>0$, \cite[1.2.12]{BigardKeimelEtAl1977} applies yielding $n \lvert g\rvert \le 1$.

Conversely, assume $n\lvert g\rvert\le 1$ for all $n\in\Npos$.
Then, taking $p=1$ and $q=n$ in the defining set we obtain
$\|g\|_\ell\le \frac1n$ for all $n$. Hence $\|g\|_\ell=0$.

For the final part of the statement, notice that, by definition, $d_\ell$ is a metric if and only if $d_\ell(g,h)=0$ implies $g=h$, i.e.\ if and only if
$\|g\|_\ell=0$ implies $g=0$.
By what we proved above, $\|g\|_\ell=0$ if and only if $n\lvert g\rvert\le 1$ for all $n$.
Thus $d_\ell$ is a metric if and only if there are no non-zero infinitesimals with respect to $1$, and by \Cref{r:unital-simplifies-Archimedean} this is equivalent to the Archimedean property for unital $\ell$-groups.
\end{proof}

\begin{lemma}\label{l:MV-metric}
Let $(G,1)$ be a unital Abelian $\ell$-group and let $A=\Gamma(G,1)$.
\begin{enumerate}
\item\label{l:MV-metric:item1} The pseudonorm $\norml{\cdot}\colon G\to\mathbb R$ restricts to a map
\begin{align*}
    \normMV{\cdot}\colon A & \longrightarrow [0,1]\\
    a & \longmapsto \normMV{a} \df \norml{a}
\end{align*}
\item\label{l:MV-metric:item2} The pseudometric $\dl(g,h)\df \norml{g-h}$ restricts to a pseudometric
\[ \dMV(a,b)\df \dl(a,b)=\norml{a-b}\qquad (a,b\in A). \]
\item\label{l:MV-metric:item3} The pseudometric $\dMV$ is a metric on $A$ if and only if the MV-algebra $A$ is semisimple.
\end{enumerate}
\end{lemma}
\begin{proof}
\eqref{l:MV-metric:item1}. Since $A=\{g\in G\mid 0\le g\le 1\}$, for every $a\in A$ we have $-1\le a\le 1$. Hence $a$ satisfies the defining condition of \eqref{eq:norml} with $p=q=1$, and so $\norml{a}\le 1$. Also, we have $\norml{a}\ge 0$ by definition of infimum. Therefore, $\norml{a}\in[0,1]$, and so the restriction $\normMV{a}\df \norml{a}$ is well-defined.

\eqref{l:MV-metric:item2}. Since $\dl$ is a pseudometric on $G$, its restriction to the subset $A\subseteq G$ is automatically a pseudometric: symmetry and the triangle inequality are inherited, and $\dMV(a,a)=\dl(a,a)=\norml{0}=0$.

\eqref{l:MV-metric:item3}. This is a consequence of \cref{lem:kernel-unit-norm,l:gamma:Archimedean}.
\end{proof}

\begin{lemma} \label{l:d=c}
    Let $A$ be an MV-algebra.
    For all $a,b \in A$,
    \[
        \dMV{(a,b)} = \normMV{\cha(a,b)}
    \]
\end{lemma}
\begin{proof}
    Working in $\Xi(A)$, we have
    \begin{align*}
        \cha(a,b) & = (a \ominus b) \lor (b \ominus a)\\
        & = ((a - b) \lor 0) \lor ((b - a) \lor 0)\\
        & = \lvert a - b\rvert \lor 0\\
        & = \lvert a - b\rvert & \text{by \eqref{eq:absolute-value-is-positive}}.
    \end{align*}
    Then,
    \begin{align*}
        \normMV{\cha(a,b)} & = \norml{\cha(a,b)} \\
        & = \norml{\lvert a - b \rvert}\\
        & = \norml{a - b} & \text{by \eqref{eq:norm-is-norm-of-absolute-value}}\\
        &= \dMV{(a,b)} &\text{by def.}&\qedhere
    \end{align*}
\end{proof}

In the rest of the paper, when we speak about convergence in a semisimple MV-algebra $A$, we will always refer to the metric $\dMV$.

\begin{remark}\label{r:sup-norm-representation}
Straightforward calculations show that when a semisimple MV-algebra is represented as an algebra of functions from a set $X$ into $[0,1]$, the norm coincides with the \emph{sup norm}: 
\[
\normMV{f} = \sup \{ f(x) \mid x \in X \}.
\]
Consequently, 
\begin{align*}
\dMV{(f,g)} = \sup \{ \lvert f(x)-g(x) \rvert \mid x \in X\}.
\end{align*}
\end{remark}

\begin{remark}
Let $(G,1)$ be a unital Abelian $\ell$-group. For $g\in G$, $p\in\Z$ and $q\in\Npos$,
we define 
\[g\preceq \frac pq \iff \lm{q}g\le \lm{p}1.\]
Notice that the above relation is well-defined (meaning that it does not depend on the pair $(p,q)$ as a representative of the rational number $\frac{p}{q}$), as the following stronger property holds: for all $p,p' \in \Z$, $q,q' \in \Npos$ and $g \in G$,
\begin{align}\label{eq:preceq-well-defined}
g \preceq \frac{p}{q} \leq \frac{p'}{q'} \qquad\text{ implies }\qquad
g \preceq \frac{p'}{q'}.
\end{align}
Indeed, $g \preceq \frac{p}{q}$ means $\lm{q}g\le \lm{p}1$, and $\frac{p}{q} \leq \frac{p'}{q'}$ means $q'p \leq qp'$.
Then, $q'qg\le q'p1 \leq qp'1$.
Since $q$ is a positive natural number, by \cite[1.2.12]{BigardKeimelEtAl1977} we can remove the two occurrences of $q$ and get $q'g \leq p'1$.

Similarly, for any MV-algebra  $A$ and any $a\in A$ we define
\[
a\preceq \frac pq \iff \lm{q}a\le \lm{p}1\text{ in }\Xi(A).
\] 
\end{remark}

Note that, by definition, for every $a$ in an MV-algebra $A$ we have
\begin{equation} \label{eq:norm-definition}
    \normMV{a} \df \inf \left\{\frac{p}{q} \in \Q\mid p,q \in \Nzero,\, q\neq 0\text{ and }a\preceq \frac{p}{q}\right\}.
\end{equation}

\begin{lemma}\label{l:leqisleq}
For every MV-algebra $A$, for every $a \in A$, for every $p \in \Nzero$ and $q \in \Npos$   
\[a \preceq \frac{p}{q}\qquad\text{ implies }\qquad\normMV{a}\leq \frac{p}{q}.\] 
In addition, if $A$ is semisimple, then the converse of the above implication holds, as well.
\end{lemma}

\begin{proof}
By \eqref{eq:norm-definition}, if $a \preceq \frac{p}{q}$ then clearly $\normMV{a}\leq\frac{p}{q}$.
Suppose now that $A$ is semisimple.
Then, by \cref{p:semisimple-algebras-of-functions}, $A$ is an algebra of functions from some set $X$ into $[0,1]$, and the norm $\normMV{\cdot}$ is the sup norm.  Thus, $\normMV{a}\leq\frac{p}{q}$ means that for every $x\in X$ we have $a(x)\preceq\frac{p}{q}$, and this is equivalent to the fact that $a \preceq \frac{p}{q}$ holds in $A$.
\end{proof}

\begin{lemma} \label{l:rewrite}
For every MV-algebra $A$, for every $a \in A$, and for every $n \in \Npos$,
\[
a \preceq \frac{1}{n} \iff \mvm{n-1}a \leq \lnot a.
\]
\end{lemma}

\begin{proof}
    In the unital $\ell$-group $\Xi(A)$, the condition $\mvm{n-1}a \leq \lnot a$ means
    \begin{equation} \label{eq:rewrite}
        (n-1)a \land 1 \leq 1 -a.
    \end{equation}
    In turn, \eqref{eq:rewrite} is equivalent to $((n-1) a \land 1) + a \leq 1$, which is equivalent to
    \[
    na \land (1 + a) \leq 1.
    \]
    We claim that this is equivalent to $na \leq 1$ (using that $a \geq 0$).
    The right-to-left implication holds trivially.
    For the left-to-right implication, since the class of Abelian $\ell$-groups is generated by $\R$ as a quasi-variety (\Cref{t:quasivariety-generated-by-linear}), it is enough to observe that in $\R$ the quasi-equation
    \[
    \forall x\, \forall u\, (x \geq 0, \, u \geq 0 ,\, nx \land (u + x) \leq u) \Rightarrow nx \leq u
    \]
    holds.
    This proves our claim.
    The condition $na \leq 1$ is equivalent to $a \preceq \frac{1}{n}$, by definition of the latter.
\end{proof}

\begin{definition}[Metrically complete $\ell$-groups, MV-algebras]\hfill
    \begin{enumerate}
        \item 
        A \emph{metrically complete $\ell$-group} is an Archimedean (necessarily Abelian) unital $\ell$-group that is complete with respect to the metric $\dl$ defined in \cref{d:ell-norm}.
        \item 
        A \emph{metrically complete MV-algebra} is a semisimple MV-algebra that is complete with respect to the metric $\dMV$ defined in \cref{l:MV-metric}.
    \end{enumerate}    
\end{definition}

\begin{notation}
We consider sequences indexed by $\Npos$ rather than $\Nzero$ only to ease the notation in the coming developments.
We write bold letters for sequences indexed by $\Npos$ and use their indexed light-faced version for elements of that sequence, so $\bm{a}$ denotes $\langle a_{n}\mid n \in \Npos\rangle$.  Also, we write $\init{a}$ for $a_{1},\dots,a_{n}$. Finally, if $\succ{a}$ and $\succ{b}$ are sequences in an MV-algebra, we write $\succ{a}\oplus\succ{b}$ for the sequence $\langle a_{1}\oplus b_{1},a_{2}\oplus b_{2},\dots\rangle$; we use a similar notation for other operations.
\end{notation}

\begin{theorem} \label{t:gamma-restricts-to-metrically-complete}
A unital Abelian $\ell$-group $G$ is metrically complete if and only if the MV-algebra $\Gamma(G,1)$ is metrically complete.
\end{theorem}
\begin{proof}
By \cite[13.2.6]{BigardKeimelEtAl1977}, if $G$ is Archimedean then it can be represented as a unital $\ell$-subgroup of $\R^X$ for some set $X$, with the unit $1$ of $G$ identified with the function $X\to\R$ constantly equal to $1$; the metric induced by the unit is therefore the uniform metric. Since we have already established that semisimple MV-algebras correspond to Archimedean unital $\ell$-groups (\cref{l:gamma:Archimedean}), we can assume that $G$ is Archimedean and hence identify $G$ with the above subgroup of $\R^X$.

 ($\Rightarrow$). Suppose  $G$ is metrically complete. Let $\succ{f}$ be a Cauchy sequence in $\Gamma(G,1)$. Then $\succ{f}$ is also a Cauchy sequence in $G$, and therefore it admits a limit $f$. Since uniform convergence implies pointwise convergence and for all $x \in X$ and $n \in \Npos$ we have $0\leq f_{n}(x) \leq 1$, it follows that $0\leq f(x) \leq 1$. Thus $f \in \Gamma(G,1)$, and $f_n$ converges to $f$ in $\Gamma(G,1)$. This shows that $\Gamma(G,1)$ is complete with respect to the uniform metric.

 ($\Leftarrow$). Suppose $\Gamma(G,1)$ is metrically complete. For every $k \in \Z$, set 
\[
\Gamma_k(G) \df \{ g \in G \mid k \leq g \leq k + 1 \}.
\] 
We have $\Gamma_0(G) = \Gamma(G,1)$.
Since $\Gamma_k(G)$ is obtained by translation of $\Gamma(G,1)$ by $k$, and since the uniform metric is translation invariant, we have that, for each $k \in \Z$, $\Gamma_k(G)$ is complete in the uniform metric. Let $\succ{f}$ be a Cauchy sequence in $G$. Since $G$ is unital, all its elements are bounded. 
By the completeness of the metric space of bounded functions from $X$ to $\R$ (with respect to the uniform metric), $\succ{f}$ converges to a certain bounded function $f \colon X \to \R$.
We show that there is $M\in \Nzero$ such that, for all $x \in X$ and $n \in \Npos$, $-M\leq f_n(x) \leq M$ (i.e., $\succ{f}$ is uniformly bounded).
Let $a \leq b \in \R$ be such that the image of $f$ is contained in $[a, b]$. 
Pick any $\varepsilon > 0$.
Since $\succ{f}$ converges to $f$, there is $n \in \Npos$ such that, for any $n' \geq n$, the image of $f_{n'}$ is contained in $[a - \varepsilon, b + \varepsilon]$.
Since $f_1, \dots, f_{n-1}$ are finitely many bounded functions, there are $c \leq d\in \R$ such that the images of $f_1, \dots, f_{n-1}$ are contained in $[c,d]$.
Pick any $M \in \Nzero$ such that $[c,d] \cup [a - \varepsilon, b + \varepsilon] \subseteq [-M, M]$; then for all $x \in X$ and $n \in \Npos$, $-M\leq f_n(x) \leq M$. Again, since uniform convergence implies pointwise convergence, for any $x \in X$ we have $-M \leq f(x) \leq M$. Since, for every $k \in \Z$, the function
\begin{align*}
    \R & \longrightarrow \R\\
    x & \longmapsto (k \lor x) \land (k + 1)
\end{align*}
is non-expansive, the sequence $\succustom{(k \lor f_n) \land (k+1)}$ is Cauchy; by completeness of $\Gamma_k(G)$, it converges to a certain $g_k \in \Gamma_k(G)$. By continuity of $\lor$ and $\land$, $g_k = (k \lor f) \land (k+1)$. For all $x \in \R$ and $m \in \Nzero$ with $0\leq x\leq m$, computation shows that
\[
x = \sum_{k=0}^{m-1} \Bigl(\big( 0 \lor (x-k) \big) \land 1\Bigr).
\]
Therefore, using the inequalities $-M\leq f \leq M$, we have
\begin{align*}
f + M &= \sum_{k=0}^{2M-1} \Bigl(\bigl( 0 \lor (f + M - k) \bigr) \land 1 \Bigr)= \sum_{k=-M}^{M-1} \Bigl(\bigl( 0 \lor (f - k) \bigr) \land 1\Bigr) \\
&=\sum_{k=-M}^{M-1} \Bigl( \bigl((k \lor f) \land (k+1)\bigr) - k \Bigl) = \sum_{k=-M}^{M-1} (g_k - k).
\end{align*}
Hence, $f = -M + \sum_{k=-M}^{M-1} (g_k - k)$.
Since $g_k \in \Gamma_k(G) \subseteq G$ and $G$ is a unital $\ell$-group, $f \in G$.
\end{proof}

\begin{definition}[$\cG$, $\cMV$]\hfill
    \begin{enumerate}
        \item We let $\cG$ denote the category of metrically complete $\ell$-groups and unital $\ell$-group homomorphisms.
        \item We let $\cMV$ denote the category of metrically complete MV-algebras and MV-homomorphisms.
    \end{enumerate}
\end{definition}

\begin{corollary} \label{c:gamma:complete}
The pair of functors $\Gamma$ and $\Xi$ of \Cref{th:gamma} restricts to an equivalence between the full subcategories $\cG$ and $\cMV$ of metrically complete $\ell$-groups and metrically complete MV-algebras.
\end{corollary}

\begin{proof}
This is an immediate consequence of \Cref{l:gamma:Archimedean} and \Cref{t:gamma-restricts-to-metrically-complete}.
\end{proof}

\section{Fast sequences} \label{s:supercauchy}

\begin{definition}[Fast sequence] \label{d:supercauchy} 
A sequence $\succ{b}$ in a metric space  $(B,d)$ is said to be \emph{fast} if for every $n \in \Npos$ we have $d(b_n,b_{n+1}) \leq \frac{1}{2^{n+1}}$.
\end{definition}

\begin{definition}[MV-fast sequence] \label{d:supercauchy-MValgebra}
A sequence $\succ{b}$ in an MV-algebra is said to be \emph{MV-fast} if for every $n \in \Npos$ we have $\cha(b_n,b_{n+1}) \preceq \frac{1}{2^{n+1}}$.
\end{definition}

Note that the condition $\cha(b_n,b_{n+1}) \preceq \frac{1}{2^{n+1}}$ appearing in \Cref{d:supercauchy-MValgebra} is not equivalent, in general, to the condition $\dMV{(b_n,b_{n+1})} \leq \frac{1}{2^{n+1}}$.
However, these two conditions are equivalent in any semisimple MV-algebra, i.e., when $\dMV$ is a metric.
The condition in \Cref{d:supercauchy-MValgebra} works better for our purposes.

\begin{lemma} \label{l:Cauchy}
Let $(B,d)$ be a metric space.
Every fast sequence in $B$ is Cauchy.  Furthermore, if $\succ{b}$ is a fast sequence converging to $b \in B$, then for every $n \in \Npos$
\[
    d(b_n,b) \leq \frac{1}{2^n}.
\]	
\end{lemma}
\begin{proof}
If $\succ{b}$ is a fast sequence in $(B,d)$, then for every $n,m \in \Npos$ such that $m\geq n$ the triangle inequality implies
\[
    d( b_n,b_m) \leq \sum_{i=n}^{m-1}d( b_i,b_{i+1}) \leq \sum_{i=n}^{m-1}\frac{1}{2^{i+1}}<\sum_{i=n +1}^{+\infty}\frac{1}{2^{i}} = \frac{1}{2^{n}}.
\]
Thus, for any $\varepsilon>0$, there is $k\in \Nzero$ such that the distance between any $b_m$  and $b_n$ with $k\leq m,n$ is smaller than $\varepsilon$. Hence, $\succ{b}$ is a Cauchy sequence.
For the second claim, suppose that the fast sequence $\succ{b}$ converges to $b \in B$, and let $n \in \Npos$.
For every $m \geq n$, 
\[
d(b_n, b) \leq d(b_n, b_m) + d(b_m, b) \leq \frac{1}{2^n} + d(b_m, b).
\]
For $m \to \infty$ we have $d(b_m, b) \to 0$;
therefore, 
\[
d(b_n, b) \leq \frac{1}{2^n}. \qedhere
\]
\end{proof}

\begin{lemma} \label{l:criterion for Cauchy}
A metric space $(B,d)$ is Cauchy complete if and only if every fast sequence converges.
\end{lemma}
\begin{proof}
By \Cref{l:Cauchy}, every fast sequence is a Cauchy sequence; thus, the left-to-right implication holds. For the converse, suppose that every fast sequence in $B$ converges. Let $\succ{b}[k]$ be a Cauchy sequence; then for every $n \in \Npos$ there is $k_n \in \Npos$ such that $d(b_{k_n},b_m) \leq \frac{1}{2^{n+1}}$ for every $m\geq k_n$. Furthermore, the sequence $\succ{k}$ can be chosen to be increasing. It follows that for every $n \in \Npos$ we have $k_{n+1} \geq k_n$; thus, $d(b_{k_n},b_{k_{n+1}}) \leq \frac{1}{2^{n+1}}$. In other words, $\succ{b_k}$ is fast and converges by the hypothesis. Since the sequence $\succ{b}$ is Cauchy and has a converging subsequence, it converges.
\end{proof}

\begin{definition}\label{d:theta}
For each $n \in \Npos$, we define the  MV-algebraic term  
\begin{align} \label{eq:def:theta}
\theta_n(x) \df x \land \Bigl(\lnot\bigl(\mvm{2^{n+1}-1} x\bigr) \Bigr).
\end{align}
\end{definition}
\begin{figure}[h!]
    \centering
    \pgfplotsset{every axis/.append style={tick label style={font=\Large}}}
\begin{tikzpicture}[scale=0.52,
declare function={
	func(\x) = (\x <1/4)*(\x)+
	and(\x >= 1/4,\x<1/3)*(1-3*\x) +
	(\x >=1/3)*(0)
	;
}
]
\begin{axis}[xmin=0,xmax=1,ymin=0,ymax=1, samples at={0,0.25,0.333333333333,1}]
\addplot[blue, ultra thick] {func(x)};
\end{axis}
\end{tikzpicture}
\begin{tikzpicture}[scale=0.52,
declare function={
	func(\x) = (\x <1/8)*(\x)+
	and(\x >= 1/8,\x<1/7)*(1-7*\x) +
	(\x >=1/7)*(0)
	;
}
]
\begin{axis}[xmin=0,xmax=1,ymin=0,ymax=1, samples at={0,0.125,0.1428571428571429,1}]
\addplot[blue, ultra thick] {func(x)};
\end{axis}
\end{tikzpicture}
\begin{tikzpicture}[scale=0.52,
declare function={
	func(\x) = (\x <1/16)*(\x)+
	and(\x >= 1/16,\x<1/15)*(1-15*\x) +
	(\x >=1/15)*(0)
	;
}
]
\begin{axis}[xmin=0,xmax=1,ymin=0,ymax=1, samples at={0,0.0625, 0.066666666666666,1}]
\addplot[blue, ultra thick] {func(x)};
\end{axis}
\end{tikzpicture}
\caption{The plots in $[0,1]$ of $\theta_n$ for $n=1,2,3$.}
\end{figure}
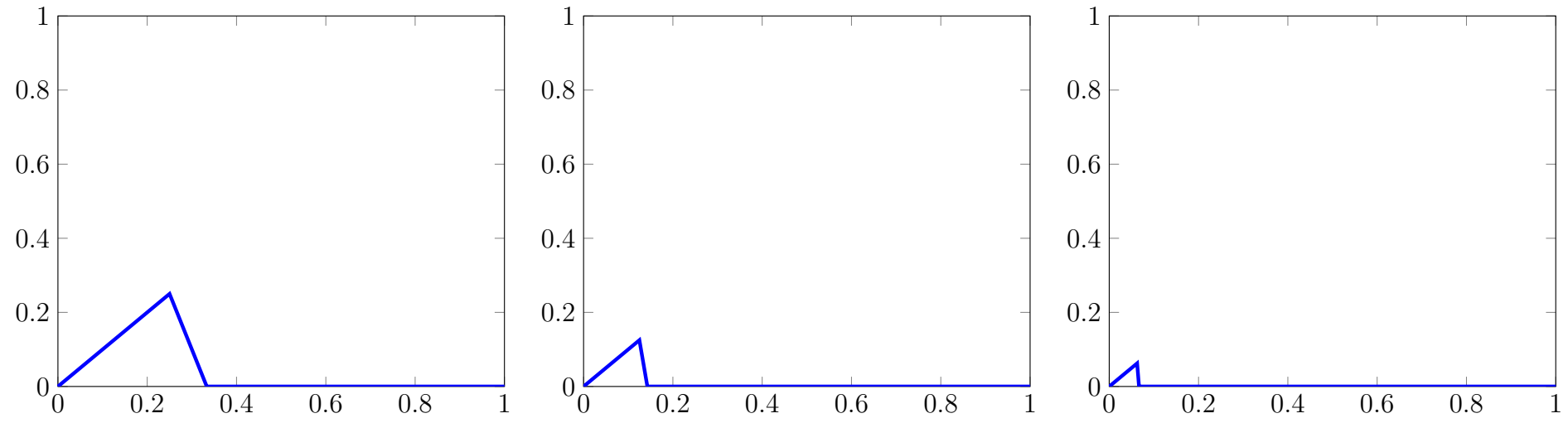

The only properties of the terms $\theta_n$ that we will use lie in the following lemma.

\begin{lemma} \label{l:bdd-norm-identity}
 Let $A$ be an MV-algebra, let $a \in A$, and let $n \in \Npos$. Then
\begin{enumerate}
    \item\label{e:little} $\theta_n(a) \preceq \frac{1}{2^{n+1}}$;
    \item\label{e:theta_little} $\theta_n(a) = a$ if and only if $a\preceq \frac{1}{2^{n+1}}$.
\end{enumerate}
\end{lemma}
\begin{proof}
\eqref{e:little}. By \cref{l:rewrite} and \Cref{l:MV-basic-properties}\eqref{l:MV-basic-properties:item-neg-antitone}, for every $x \in A$
\begin{align}\label{eq:rewrite2}
x\preceq \frac{1}{2^{n+1}}\iff\mvm{2^{n+1}-1}x\leq \neg x \iff x\leq  \neg\mvm{2^{n+1}-1}x.
\end{align}
Therefore, the condition in \eqref{e:little} is equivalent to an equation, and, by \cref{l:MV-subd-lin-ord}, it is enough to prove it in all linearly ordered MV-algebras. 
If the order is linear, we have two cases: either $a\leq  \neg\mvm{2^{n+1}-1}a$, and thus $\theta_{n}(a) = a\preceq \frac{1}{2^{n+1}}$ and so we are done; or $\neg\mvm{2^{n+1}-1}a\leq a$, and thus
$\mvm{2^{n+1}-1}(\neg\mvm{2^{n+1}-1}a)\leq \mvm{2^{n+1}-1}a$. In the latter case, we have $\theta_{n}(a) = \neg\mvm{2^{n+1}-1}a$, and an application of \cref{l:rewrite} to the element $\neg\mvm{2^{n+1}-1}{a}$ gives
\[\theta_n(a)=\neg\mvm{2^{n+1}-1}{a}\preceq \frac{1}{2^{n+1}}.\]

\eqref{e:theta_little}.
If $\theta_n(a) = a$, then $a = \theta_n(a) \preceq \frac{1}{2^{n+1}}$ by the previous item.
Conversely, if $a\preceq \frac{1}{2^{n+1}}$, then by \eqref{eq:rewrite2} we have $a\leq\neg\mvm{2^{n+1}-1}a$, and so $a=a \land (\lnot\mvm{2^{n+1}-1}a) = \theta_n(a)$ by definition of $\theta_n$ (\cref{eq:def:theta}).
\end{proof}

We now use these $\theta_n$ to define a family of terms that transforms any sequence into an MV-fast sequence. Furthermore, if the sequence is already MV-fast, then the terms do not alter it.
\begin{definition}
The family of MV-algebraic terms $\{\rho_{n}\mid n \in \Npos\}$, where each $\rho_{n}$ has arity $n$, is defined recursively as follows:
\begin{align*}
\rho_1(x_1)&\df x_1,\\
\rho_{n+1}(\init{x}[n+1])&\df \bigl(\rho_{n}(\init{x}) \oplus \theta_{n}(x_{n+1}\ominus x_n) \bigr) \ominus \theta_{n}(x_{n}\ominus x_{n+1}).
\end{align*}
We write $\asucc{\rho(x)}$ for the sequence $\succustom{\rho_{n}(\init{x})}$.
\end{definition}

\begin{lemma} \label{l:in V dist rho_n rho_m}
Let $A$ be an MV-algebra and $\succ{a}$ be a sequence in $A$. The following properties hold:
\begin{enumerate}
\item\label{i:rho_super_C} the sequence $\succ{\rho(a)}$ is MV-fast in $A$;
\item\label{i:rho_proj} if $\succ{a}$ is MV-fast, then, for all $n \in \Npos$,  $\rho_n(\init{a}) = a_n$.
\end{enumerate}
\end{lemma}
\begin{proof}
\eqref{i:rho_super_C}.
From \Cref{l:MV-basic-properties}\eqref{l:MV-basic-properties:item-monotone} it immediately follows that, for any $n \in \Npos$,
\[
\rho_{n}(\init{a}) \ominus\theta_{n}(a_{n}\ominus a_{n+1}) \leq\rho_{n+1}(\init{a}[n+1]) \leq \rho_{n}(\init{a}) \oplus\theta_n( a_{n+1}\ominus a_n).
\]	
The displayed inequalities entail
\begin{align*}
&\rho_{n+1}(\init{a}[n+1])\ominus \rho_n(\init{a})\\
& \le \big(\rho_n(\init{a})\oplus \theta_n(a_{n+1}\ominus a_n)\big)\ominus \rho_n(\init{a}) & \text{by \Cref{l:MV-basic-properties}\eqref{l:MV-basic-properties:item-monotone}}\\
&\le \theta_n(a_{n+1}\ominus a_n) & \text{because }
(x\oplus y)\ominus x\leq y.
\end{align*}
Similarly,
\begin{align*}
    &\rho_n(\init{a})\ominus \rho_{n+1}(\init{a}[n+1]) \\
    &\le \rho_n(\init{a}) \ominus (\rho_{n}(\init{a}) \ominus\theta_{n}(a_{n}\ominus a_{n+1})) & \text{by \Cref{l:MV-basic-properties}\eqref{l:MV-basic-properties:item-monotone}}\\
    & \leq \theta_n(a_n\ominus a_{n+1}) & \text{because } x\ominus(x\ominus y)\leq y.
\end{align*}
Since \(u\preceq r\) and \(v\preceq r\) together imply \(u\vee v\preceq r\), we have
\[
\cha(\rho_{n}(\init{a}),\rho_{n+1}(\init{a}[n+1]))
\le \theta_n(a_{n+1}\ominus a_n) \vee \theta_n(a_n\ominus a_{n+1})
\preceq \frac{1}{2^{n+1}},
\]
where the last inequality follows from \Cref{l:bdd-norm-identity}\eqref{e:little}.
Since $a \leq b \preceq q$ implies $a \preceq q$, we obtain $\cha(\rho_{n}(\init{a}),\rho_{n+1}(\init{a}[n+1])) \preceq \frac{1}{2^{n+1}}$.

\eqref{i:rho_proj}. Let $\succ{a}$ be MV-fast. We prove,  by induction on $n \in \Npos$, that $\rho_n(\init{a}) = a_n$. By definition, $\rho_1(a_1) = a_1$. Suppose that, for a fixed $n \in \Npos$, $\rho_n(\init{a}) = a_n$. Since
\[
\cha(a_n,a_{n+1}) \preceq \frac{1}{2^{n+1}}
\]
and 
\[
\cha(a_n,a_{n+1}) = (a_{n+1}\ominus a_n) \lor (a_n\ominus a_{n+1}),
\]
we have 
\[a_{n+1}\ominus a_n \preceq \frac{1}{2^{n+1}}\text{ and }a_n\ominus a_{n+1} \preceq\frac{1}{2^{n+1}}.\] 
By \Cref{l:bdd-norm-identity}\eqref{e:theta_little}, $\theta_n(a_{n+1}\ominus a_n) = a_{n+1}\ominus a_n$ and $\theta_n(a_{n}\ominus a_{n+1}) = a_{n}\ominus a_{n+1}$. Then
\begin{align*}
\rho_{n+1}(\init{a}[n+1])& = (\rho_{n}(\init{a}) \oplus \theta_{n}(a_{n+1}\ominus a_n)) \ominus \theta_{n}(a_{n}\ominus a_{n+1}) & \text{by definition}\\
& = (a_n\oplus  (a_{n+1}\ominus a_n)) \ominus (a_{n}\ominus a_{n+1})& \text{induct.\ hypothesis}\\
&=a_{n+1} & \text{by \Cref{l:MV-basic-properties}\eqref{l:MV-basic-properties:item-cancellation-identity}}. &\qedhere
\end{align*}
\end{proof}

\begin{lemma} \label{l:fast-iff}
Let $A$ be an MV-algebra. A sequence $\succ{a}$ in $A$ is MV-fast if and only if there is $\succ{b}$ in  $A$ such that, for every $n \in \Npos$, $a_n = \rho_n(\init{b})$.
\end{lemma}
\begin{proof}
($\Rightarrow$). If $\succ{a}$ in $A$ is MV-fast, then set $b_n\df a_n$ for each $n \in \Npos$. By \Cref{l:in V dist rho_n rho_m}\eqref{i:rho_proj}, for every $n \in \Npos$ we have $a_n = \rho_n(\init{b})$.

($\Leftarrow$). This is \Cref{l:in V dist rho_n rho_m}\eqref{i:rho_super_C}.
\end{proof}

Let $A$ be a metrically complete MV-algebra and let $\succ{a}$ be a sequence in $A$. By \Cref{l:in V dist rho_n rho_m}, $\succ{\rho(a)}$ is a fast sequence, and so it is also Cauchy by \Cref{l:Cauchy}. Thus, it admits a unique limit, which we denote by 
\[
\limn\rho_n(\init{a}).
\]
\begin{remark} \label{l:gamma is lim}
If $\succ{a}$ is a fast sequence in a metrically complete MV-algebra $A$, then
\[
\limn\rho_n(\init{a}) = \limn a_n.
\]
Indeed, since $\succ{a}$ is fast, for all $n \in \Npos$ we have $\rho_n(\init{a}) = a_n$, by \Cref{l:in V dist rho_n rho_m}\eqref{i:rho_proj}.
\end{remark}

\begin{proposition} \label{p:then CMV}
Let $A$ be a metrically complete MV-algebra. 
\begin{enumerate}
\item\label{i:gamma1} For $a \in A$, $\limn\rho_n(a,\dots,a) = a$.
\item\label{i:gamma2} For $a \in A$, $\limn\rho_n(\theta_1(a),\dots, \theta_n(a)) = 0$.
\item\label{i:gamma3} For every sequence $\succ{a}$ in $A$ and $m \in \Npos$, \[\cha\mleft(\limn\rho_n(\init{a}),\rho_{m}(\init{a}[m])\mright) \preceq \frac{1}{2^m}.\]
\end{enumerate}
\end{proposition}
\begin{proof}
The sequence constantly equal to $a$ is clearly MV-fast. Thus, by \Cref{l:gamma is lim} we have $\limn\rho_n(a,\dots,a) = \limn a=a$.  This proves \cref{i:gamma1}.

To prove \cref{i:gamma2} we calculate:
\begin{align*}
\cha(\theta_n(a), \theta_{n+1}(a)) & = (\theta_n(a) \ominus \theta_{n+1}(a)) \lor ( \theta_{n+1}(a) \ominus \theta_n(a))&\text{by definition of $\cha$}\\
&\leq \theta_n(a) \lor \theta_{n+1}(a)&\text{since $x \ominus y \leq x$}\\
&\preceq \frac{1}{2^{n+1}}& \text{by \Cref{l:bdd-norm-identity}\eqref{e:little}}.&
\end{align*}
So, $\succustom{\theta_n(a)}$ is MV-fast. Hence, we can apply \Cref{l:gamma is lim} to obtain 
\[\limn\rho_n(\theta_1(a),\dots,\theta_n(a)) = \limn\theta_n(a).\] 
To conclude the proof of \cref{i:gamma2}, notice that by \Cref{l:bdd-norm-identity}\eqref{e:little} and \cref{l:leqisleq} we have $\limn\normMV{\theta_n(a)} = 0$. Therefore, $d_{\mathrm{MV}}(\theta_n(a),0)=\|\theta_n(a)\|_{\mathrm{MV}}\to0$. Hence, $\limn\theta_n(a) = 0$.

Finally, \cref{i:gamma3} is a direct consequence of \Cref{l:Cauchy}, as it implies that,  for $m \in \Npos$, 
\begin{align} \label{eq:p:then CMV}
\left\lVert\cha\mleft(\limn\rho_n(\init{a}),\rho_{m}(\init{a}[m])\mright)\right\rVert_{\mathrm{MV}}
=\dMV\mleft(\limn\rho_n(\init{a}),\rho_{m}(\init{a}[m]) \mright) \leq \frac{1}{2^m},
\end{align} 
where the equality follows from \cref{l:d=c}.
So, it is enough to apply the equivalence under semisimplicity in \cref{l:leqisleq} to obtain $\cha\mleft(\limn\rho_n(\init{a}),\rho_{m}(\init{a}[m])\mright) \preceq \frac{1}{2^m}$.
\end{proof}

\section{The infinitary variety \texorpdfstring{$\CMV$}{CMV}, and its equivalence with metrically complete MV-algebras} \label{s:CMV}
\begin{definition}[The infinitary variety $\CMV$] \label{d:CMV}
We denote by $\CMV$ the (infinitary) variety obtained from the variety of MV-algebras by adding an operation $\gamma$ of countably infinite arity, subject to the following additional axioms.
\begin{enumerate}[label = (A\arabic*), ref = A\arabic*]
\item\label{d:CMV:A1} $\gamma(x,x,\dots) = x$.
\item\label{d:CMV:A2} $\gamma(\theta_1(x), \theta_2(x),\dots) = 0$.
\item\label{d:CMV:A3} For each $n \in \Npos$, $\cha\mleft(\gamma(\asucc{x}),\rho_{n}(\init{x}) \mright) \preceq \frac{1}{2^n}$.
\end{enumerate} 
An object of $\CMV$ is called a \emph{$\CMV$-algebra}.
\end{definition}
This is indeed a variety because, by \Cref{l:rewrite}, the axiom \eqref{d:CMV:A3} is equivalent to a family of equations:
\[\mvm{2^n - 1}\cha\mleft(\gamma(\asucc{x}),\rho_{n}(\init{x}) \mright) \leq \neg\cha\mleft(\gamma(\asucc{x}),\rho_{n}(\init{x}) \mright),\qquad n \in \Npos.\]

\begin{proposition} \label{p:CMV arch}
Every MV-algebra $A$ that can be equipped with a function
\[
\gamma \colon A^\Npos\to A
\]
satisfying \eqref{d:CMV:A1} and \eqref{d:CMV:A2} is semisimple. In particular,
every $\CMV$-algebra is a semisimple MV-algebra.
\end{proposition}
\begin{proof}
We use the following characterisation of semisimple MV-algebras (see \cite[Prop.~3.6.4]{MundiciAlg}): an MV-algebra $A$ is semisimple if and only if, for every $a \in A$,
\begin{align} \label{eq:infinitesimal}
\text{if }a\leq \neg\mleft(\mvm{n} a\mright) \text{ for all }n \in \Npos\text{ then }a=0. 
\end{align} 
Suppose that $a \in A$ satisfies the hypothesis of \cref{eq:infinitesimal}.  It follows that $a\leq \lnot\mleft(\mvm{2^{n+1}-1}a\mright)$  for all $n \in \Npos$; hence $\theta_n(a) = a \land \lnot\mleft(\mvm{2^{n+1}-1}a\mright) = a$. Therefore, using \eqref{d:CMV:A1} and \eqref{d:CMV:A2},
\[
a = \gamma(a,a,\dots) = \gamma\bigl(\theta_1(a), \theta_2(a),\dots\bigr) = 0. \qedhere
\]
\end{proof}

\begin{remark} \label{r:eqA3}
If axioms \eqref{d:CMV:A1} and \eqref{d:CMV:A2} hold in $A$, then $A$ is semisimple by \cref{p:CMV arch}, and, by \cref{l:leqisleq}, the condition \eqref{d:CMV:A3} is equivalent to
\begin{enumerate}[label = (A3$^\prime$), ref = A3$^\prime$]
\item\label{d:CMV:A3b} For each $n \in \Npos$, $\dMV(\gamma(\asucc{x}),\rho_{n}(\init{x})) \leq \frac{1}{2^n}$.
\end{enumerate}
Thus, $\CMV$-algebras can be equivalently defined as MV-algebras equipped with an operation $\gamma$ of countably infinite arity satisfying \eqref{d:CMV:A1}, \eqref{d:CMV:A2} and \eqref{d:CMV:A3b}.
As a consequence, for any $\CMV$-algebra $A$ and for any  sequence $\succ{a}$ in $A$,
\[
\gamma(\asucc{a}) = \limn\rho_n(\init{a}).
\]
However, \eqref{d:CMV:A3b} is clearly not an equation in the language of MV-algebras.
\end{remark}

\begin{proposition} \label{p:CMV Cauchy complet}
Any $\CMV$-algebra is a metrically complete MV-algebra.	
\end{proposition}
\begin{proof}Let $A$ be a $\CMV$-algebra. By \Cref{p:CMV arch},  $A$ is semisimple. To show that $A$ is metrically complete, by \Cref{l:criterion for Cauchy} it is enough to prove that every fast sequence $\succ{a}$ in $A$ converges. By  \Cref{l:in V dist rho_n rho_m}, for every $n \in \Npos$, $\rho_n(\init{a}) = a_n$. By \Cref{r:eqA3}, for each $n \in \Npos$,
\[
\dMV(\gamma(\asucc{a}),a_n) = \dMV(\gamma(\asucc{a}),\rho_{n}(\init{a})) \leq \frac{1}{2^n}.
\]
Hence, $\succ{a}$ converges to $\gamma(\asucc{a})$.
\end{proof}

\begin{lemma} \label{l:homo preserves lim}
Any homomorphism $f \colon A\to B$ of MV-algebras is a non-expansive map, i.e.,
\[\dMV{(f(a_{1}),f(a_{2}))}\leq \dMV{(a_{1},a_{2})}.\]  
Thus, if $A$ and $B$ are semisimple MV-algebras, then $f$ is continuous with respect to the metric topologies of $A$ and $B$.
\end{lemma}
\begin{proof}
Let $\widehat f\df \Xi(f)$. Fix $a_1,a_2\in A$. By definition of $\dMV$ as the restriction of the $\ell$-group pseudometric $\dl$, and since $\widehat f$ is a group
homomorphism, we have
\[
\dMV\bigl(f(a_1),f(a_2)\bigr)
= \dl\bigl(\widehat f(a_1),\widehat f(a_2)\bigr)
= \norml{\widehat f(a_1)-\widehat f(a_2)}
= \norml{\widehat f(a_1-a_2)}.
\]

Now let $\frac{p}{q}$ be any rational with $p\in\Nzero$ and $q\in\Npos$ such that
\[
-p1 \leq q(a_1-a_2)\leq p1
\qquad\text{in }\Xi(A).
\]
Applying $\widehat f$ and using that $\widehat f$ is order-preserving and preserves the
unit, we obtain
\[
-p1 = \widehat f(-p1)\leq \widehat f\bigl(q(a_1-a_2)\bigr)
= q\widehat f(a_1-a_2)
\leq \widehat f(p1)=p1
\qquad\text{in }\Xi(B).
\]
Hence, every $\frac{p}{q}$ that witnesses the defining condition for $\norml{a_1-a_2}$
also witnesses the defining condition for $\norml{\widehat f(a_1-a_2)}$. Taking infima
over all such rationals yields
\[
\norml{\widehat f(a_1-a_2)} \leq \norml{a_1-a_2},
\]
and therefore
\[
\dMV\bigl(f(a_1),f(a_2)\bigr)\leq \dMV(a_1,a_2),
\]
as required.

Finally, if $A$ and $B$ are semisimple, then $\dMV$ is a metric on both, and any
non-expansive map between metric spaces is continuous. 
\end{proof}

\begin{proposition} \label{p:full}
Let $A$ and $B$ be $\CMV$-algebras, and let $f \colon A\to B$ be a homomorphism of MV-algebras. Then $f$ preserves $\gamma$.
\end{proposition}
\begin{proof}
Let $\succ{a}$ be a sequence of elements of $A$. Then
\begin{align*}
f(\gamma(\asucc{a}))&=f\mleft(\limn\rho_n(\init{a}) \mright) &\text{by \Cref{r:eqA3}}\\
& = \limn f(\rho_n(\init{a}))  &\text{by \Cref{l:homo preserves lim}}\\
& = \limn\rho_n(f(a_1),\dots,f(a_n)) &\text{because $f$ is a homomorphism}\\
& = \gamma(f(a_1),f(a_2),\dots) &\text{by \Cref{r:eqA3}}.& \qedhere
\end{align*}
\end{proof}

We close this section by summarising the main properties of the forgetful functor from $ \CMV$ into $\MV$.
\begin{theorem} \label{t:forget-functor}
Let $\U \colon \CMV\to \MV$ be the functor that ``forgets'' the operation $\gamma$.
\begin{enumerate}
\item\label{t:forget-functor:item1} The functor $\U$ is full and faithful.
\item\label{t:forget-functor:item2} The functor $\U$ is injective on objects: $\U(A) = \U(B)$ implies $A=B$.
\item\label{t:forget-functor:item3} The image on objects of $\U$ is the class of metrically complete MV-algebras.
\end{enumerate}
Therefore, $\U$ restricts to an isomorphism of categories between $\CMV$ and $\cMV$, the latter being the full subcategory of $\MV$ given by metrically complete MV-algebras.
\end{theorem}
\begin{proof}
\eqref{t:forget-functor:item1}. This holds because $\U$ is clearly faithful and full by \Cref{p:full}.

\eqref{t:forget-functor:item2}. Let $\U(A) = \U(B)$.  Since $\U$ is full, there is a $\CMV$-morphism $\varphi \colon A\to B$ such that $\U(\varphi) = \mathtt{id} \colon \U(A) \to \U(B)$. By the definition of $\U$, such a $\varphi$ must be the identity function; thus $A=B$. 

\eqref{t:forget-functor:item3}. Suppose that $A$ is semisimple and Cauchy complete. Then, by \Cref{p:then CMV}, $A$ admits a structure of $\CMV$-algebra that extends its structure of MV-algebra. Vice versa, if there is $\tilde{A}\in \CMV$ such that $\U(\tilde{A}) = A$, then, by \Cref{p:CMV arch,p:CMV Cauchy complet}, $A$ is metrically complete. 
\end{proof}

\begin{theorem} \label{t:cG-equiv-CMV}
The categories $\cG$ of metrically complete $\ell$-groups and $\cMV$ of metrically complete MV-algebras are equivalent to the variety of (infinitary) algebras $\CMV$.
\end{theorem}
\begin{proof}
By \Cref{c:gamma:complete}, Mundici's functor $\Gamma$ restricts to an equivalence of categories
\[
\Gamma\colon \cG \simeq \cMV.
\]

Define a functor
\[
\F\colon \cMV \longrightarrow \CMV
\]
as follows. 

On objects: if $A\in\cMV$, define an infinitary operation
$\gamma_A\colon A^{\Npos}\to A$ by
\[
{\gamma}_A(\asucc{x}) \df \lim_{n\to +\infty}\rho_n(\init{x}),
\]
where $\rho_n$ are the MV-terms introduced in \Cref{s:supercauchy}. The sequence
$\succ{\rho(x)}$ is MV-fast by \Cref{l:in V dist rho_n rho_m}. Hence, it is Cauchy by \Cref{l:Cauchy} and its limit exists because $A$ is complete. By \Cref{p:then CMV},
the operation $\gamma_A$ satisfies the axioms \eqref{d:CMV:A1}--\eqref{d:CMV:A3} of \Cref{d:CMV};
therefore, $(A,\gamma_A)$ is a $\CMV$-algebra, and we set $\F(A)\df (A,\gamma_A)$.

On morphisms: if $f\colon A\to B$ is an MV-homomorphism between objects of $\cMV$, we set
$\F(f)\df f$. This is well-defined because $f$ is non-expansive
(\Cref{l:homo preserves lim}), hence preserves limits of Cauchy sequences in complete metric spaces,
and therefore preserves $\gamma$ (cf.~\Cref{p:full}). Thus $\F(f)$ is a $\CMV$-homomorphism.

Let $\U\colon \CMV\to \MV$ be the forgetful functor. By construction,
\[
\U\circ \F = \mathrm{Id}_{\cMV}.
\]
Conversely, if $A\in\CMV$, then by \Cref{r:eqA3} the $\gamma$-operation of $A$ is uniquely determined
by the MV-reduct $\U(A)$ via the formula
\[
\gamma(\asucc{x})=\lim_{n\to+\infty}\rho_n(\init{x}).
\]
Thus, $\F(\U(A))=A$ (and similarly on morphisms).
Hence $\F$ is a quasi-inverse of
$\U\colon\CMV\to\cMV$, and therefore
\[
\F\colon \cMV \cong \CMV.
\]

Composing the equivalence with the isomorphism, we obtain the equivalence \({\F}\circ {\Gamma} \colon \cG \simeq \CMV\), as required.
\end{proof}

\section{A finite equational axiomatisation of \texorpdfstring{$\CMV$}{CMV}} \label{s:finite-axiom}
The aim of this section is to show that $\CMV$ admits a finite equational axiomatisation (in the infinitary language considered in the previous section).
\begin{notation}	
Given a sequence $\succ{a}\df (a_{1},a_{2},\dots)$, we write $\shift{\succ{a}}$ for the sequence $(a_{2},a_{3},\dots)$.
\end{notation}

\begin{lemma} \label{l:oplus-and-odot-are-fast}
Let $A$ be a semisimple MV-algebra and $\succ{a}$ be a converging sequence in $A$.  The following conditions hold:
\begin{enumerate}
\item\label{l:oplus-and-odot-limits:item1} $\lim(\succ{a}\oplus \succ{a}) = \lim\succ{a}\oplus \lim\succ{a}$;
\item\label{l:oplus-and-odot-limits:item3} $\lim (\succ{a}\odot\succ{a}) = \lim\succ{a}\odot \lim\succ{a}$.
\end{enumerate}  
\end{lemma}
\begin{proof}
It is enough to show that if $\dMV(a,b) \leq \frac{1}{2^{n+1}}$ then $\dMV(a\oplus a,b\oplus b) \leq \frac{1}{2^n}$ and $\dMV(a\odot a,b\odot b) \leq \frac{1}{2^n}$.
Since $A$ is semisimple, by \cref{l:leqisleq} the claim reduces to checking that 
$\cha(a\oplus a,b\oplus b) \preceq \frac{1}{2^n}$ and $\cha(a\odot a,b\odot b) \preceq \frac{1}{2^n}$. In turn, by \cref{l:rewrite}, these two conditions are equivalent to two equations. Thus, the claim in the lemma reduces to checking that two quasi-equations hold.  Since quasi-equations are preserved by products and subalgebras, by \Cref{p:semisimple-algebras-of-functions} it is enough to check that they hold in $[0,1]$. This can be easily seen using the following inequalities:
\[
\lvert (a \oplus a) - (b \oplus b)\rvert = \big\lvert \min\{2a,1\}-\min\{2b,1\}\big\rvert \le 2 \lvert a-b \rvert,
\]
and
\[
\lvert (a \odot a) - (b \odot b)\rvert = \big\lvert \max\{0,2a-1\}-\max\{0,2b-1\}\big\rvert \le 2 \lvert a-b \rvert. \qedhere
\]
\end{proof}

\begin{lemma} \label{l:basic-properties-norm}
For every MV-algebra $A$ and $a,b \in A$,
\[
    \cha(\mvm{2}a,\mvm{2}b)\leq \mvm{2}\cha(a,b)
    \quad\text{and}\quad
    \cha(a\odot a,b\odot b)\leq \mvm{2}\cha(a,b).
\]
\end{lemma}
\begin{proof}
Since the order of an MV-algebra is equationally definable, both inequalities are
equations in the language of MV-algebras. Hence, by \Cref{t:dinola}, it is enough
to check them in the standard MV-algebra $[0,1]$.

Let $r,s\in[0,1]$. There,
\[
    \cha(r,s)=\lvert r-s \rvert,
    \qquad
    \mvm{2}r=\min\{2r,1\},
    \qquad
    r\odot r=\max\{2r-1,0\}.
\]
The maps $t\mapsto \min\{2t,1\}$ and $t\mapsto \max\{2t-1,0\}$ are
$2$-Lipschitz. Therefore
\[
    \cha(\mvm{2}r,\mvm{2}s)\leq 2 \lvert r-s\rvert,
    \qquad
    \cha(r\odot r,s\odot s)\leq 2\lvert r-s \rvert.
\]
Since both left-hand sides are bounded by $1$, we get
\[
    \cha(\mvm{2}r,\mvm{2}s)\leq \min\{2 \lvert r-s\rvert,1\}
    =\mvm{2}\cha(r,s),
\]
and similarly
\[
    \cha(r\odot r,s\odot s)\leq \min\{2 \lvert r-s \rvert,1\}
    =\mvm{2}\cha(r,s).
\]
Thus both inequalities hold in $[0,1]$, and hence in every MV-algebra.
\end{proof}

\begin{lemma}\label{l:a-fast-implies-a+a-fast}
If $\succ{a}$ is an MV-fast sequence in an MV-algebra, then  $\shift{\succ{a}}\oplus \shift{\succ{a}}$ and $\shift{\succ{a}}\odot \shift{\succ{a}}$ are MV-fast.
\end{lemma}
\begin{proof}
Since $\succ{a}$ is MV-fast, $\cha(a_{n+1},a_{n+2})\preceq \frac{1}{2^{n+2}}$, for every $n \in \Npos$. Hence, by \cref{l:basic-properties-norm}, 
\[
\cha(a_{n+1}\oplus a_{n+1},a_{n+2}\oplus a_{n+2}) \leq \mvm{2}\cha(a_{n+1},a_{n+2})\preceq  \frac{1}{2^{n+1}}.
\]
Since $a_{n+1}\oplus a_{n+1}$ is the $n^{\rm th}$ term of the sequence ${\shift{\succ{a}}}\oplus {\shift{\succ{a}}}$, the latter is MV-fast.  Similar considerations can be made for $\shift{\succ{a}}\odot\shift{\succ{a}}$.
\end{proof}

\begin{proposition} \label{p:soundness finite}
For any $\CMV$-algebra $A$, the following properties hold:
\begin{enumerate}[label = (B\arabic*), ref = B\arabic*]
\item\label{s:finite-axiom:item1} $\forall x\ \ \gamma(x,x,x,\dots) = x$;
\item\label{s:finite-axiom:item2} $\forall x \ \ \gamma(\theta_1(x), \theta_2(x),\theta_3(x), \dots) = 0$;
\item\label{s:finite-axiom:item3} $\forall \succ{x}\ \ \cha(\gamma(\asucc{x}), x_1) \odot \cha(\gamma(\asucc{x}), x_1) = 0$;
\item\label{s:finite-axiom:item4} $\forall\succ{x}\ \  \gamma\bigl(\asucc{\rho(x)}\bigr) \oplus\gamma\bigl(\asucc{\rho(x)}\bigr) = \gamma\bigl(\shift{\asucc{\rho(x)}}\oplus\shift{\asucc{\rho(x)}}\bigr)$;
\item\label{s:finite-axiom:item5} $\forall \succ{x}\ \  \gamma\bigl(\asucc{\rho(x)}\bigr) \odot\gamma\bigl(\asucc{\rho(x)}\bigr) = \gamma\bigl(\shift{\asucc{\rho(x)}}\odot\shift{\asucc{\rho(x)}}\bigr)$;
\item\label{s:finite-axiom:item6} $\forall \succ{x}\ \ \gamma(\asucc{x}) = \gamma\bigl(\asucc{\rho(x)}\bigr)$.
\end{enumerate}
\end{proposition}
\begin{proof}
The conditions \eqref{s:finite-axiom:item1} and \eqref{s:finite-axiom:item2} are simply \eqref{d:CMV:A1} and \eqref{d:CMV:A2} in \Cref{d:CMV}, and \eqref{s:finite-axiom:item3} is an equivalent formulation of case $n=1$ of \eqref{d:CMV:A3} in \Cref{d:CMV}; hence, \eqref{s:finite-axiom:item1}, \eqref{s:finite-axiom:item2}, and \eqref{s:finite-axiom:item3} hold in any $\CMV$-algebra.

We now prove \eqref{s:finite-axiom:item4}. Let $\succ{a}$ be a sequence in $A$.  By \Cref{l:in V dist rho_n rho_m}, the sequence $\succ{\rho(a)}$ is fast; hence, by \Cref{l:a-fast-implies-a+a-fast}, the sequence $\shift{\succ{\rho(a)}}\oplus \shift{\succ{\rho(a)}}$ is fast, too. Hence, 
\begin{align*}	
\gamma\bigl(\shift{\succ{\rho(a)}}\oplus\shift{\succ{\rho(a)}}\bigr)& = \lim\bigl(\shift{\succ{\rho(a)}}\oplus \shift{\succ{\rho(a)}}\bigr) &\text{ by  \Cref{l:gamma is lim}}\\
& = \lim(\succ{\rho(a)}\oplus \succ{\rho(a)})\\
& = \lim \succ{\rho(a)}\oplus \lim \succ{\rho(a)}& \text{by \Cref{l:oplus-and-odot-are-fast}\eqref{l:oplus-and-odot-limits:item1}}\\
& = \gamma(\succ{\rho(a)}) \oplus\gamma(\succ{\rho(a)})&\text{ by  \Cref{l:gamma is lim}.}
\end{align*}

\eqref{s:finite-axiom:item5} can be proved similarly to  \eqref{s:finite-axiom:item4}. 
Finally, to prove \eqref{s:finite-axiom:item6}, recall that, on the one hand, by \Cref{r:eqA3}, in any $\CMV$-algebra we have
\[
\gamma(\asucc{x}) = \limn\rho_n(\init{x}).
\]
On the other hand, by \Cref{l:in V dist rho_n rho_m}, $\succ{\rho(x)}$ is a fast sequence. Hence, by  \Cref{l:gamma is lim}, 
\[
\gamma(\asucc{\rho(x)}) = \limn\rho_n(\init{x}).
\]
Therefore, 
\begin{equation*}
\gamma(\asucc{x}) = \limn\rho_n(\init{x}) = \gamma(\asucc{\rho(x)}). \qedhere
\end{equation*}
\end{proof}

\begin{remark}
\eqref{s:finite-axiom:item3} is equivalent to $\cha(\gamma(\asucc{x}), x_1) \preceq \frac{1}{2}$.
By \Cref{p:CMV arch},  \eqref{s:finite-axiom:item1} and \eqref{s:finite-axiom:item2} imply semisimplicity, and thus \eqref{s:finite-axiom:item3} is also equivalent to $\dMV(\gamma(\asucc{x}), x_1) \leq \frac{1}{2}$.

Furthermore, by  \Cref{l:fast-iff}, \eqref{s:finite-axiom:item4} and \eqref{s:finite-axiom:item5} are equivalent to the conditions \eqref{i:asucc-oplus} and \eqref{i:asucc-odot} below, respectively. 
\begin{enumerate}[label = (B\arabic*$^\prime$), ref = B\arabic*$^\prime$, start = 4]
\item \label{i:asucc-oplus}
For every MV-fast sequence $ \succ{x}$,  $\gamma(\asucc{x}) \oplus\gamma(\asucc{x}) = \gamma(\shift{\asucc{x}}\oplus\shift{\asucc{x}})$.
\item \label{i:asucc-odot}
For every MV-fast sequence $ \succ{x}$, $\gamma(\asucc{x}) \odot \gamma(\asucc{x}) = \gamma(\shift{\asucc{x}}\odot\shift{\asucc{x}})$.
\end{enumerate}
\end{remark}

We conclude this section by proving that every MV-algebra endowed with an additional operation $\gamma$ satisfying \eqref{s:finite-axiom:item1}--\eqref{s:finite-axiom:item6} is in $\CMV$.  This is the content of  \Cref{t:B1-B6-imply-CMV} below, for whose proof we need a number of preliminary results.
\begin{notation}	
We will need to consider iterated combinations of the operations $x\oplus x$ and $x\odot x$.  So, for $n \in \Nzero$ let us define the set 
\[T_{n} \coloneqq\{\tau_{n}(\tau_{n-1}(\dots(\tau_1(x)) \dots)) \mid \tau_{i}(x) \text{ is either }x\oplus x\text{ or }x\odot x\}.\]
\end{notation}	

\begin{figure}
    \centering
    \newcommand{\piecewisepanel}[1]{%
    \begin{subfigure}[t]{0.24\linewidth}
        \centering
        \begin{tikzpicture}[scale=0.52]
            \begin{axis}[
                scale only axis,
                width=5.15cm, height=4.25cm,
                xmin=0, xmax=1, ymin=0, ymax=1,
                xtick={0,0.25,0.5,0.75,1},
                ytick={0,1},
                tick label style={font=\Large},
                samples at={0,0.25,0.5,0.75,1}
            ]
                \addplot[blue, ultra thick, no marks] {min(1,max(0,4*x-#1))};
            \end{axis}
        \end{tikzpicture}
        \caption{$k = #1$}
    \end{subfigure}%
}
\piecewisepanel{0}\hfill
\piecewisepanel{1}\hfill
\piecewisepanel{2}\hfill
\piecewisepanel{3}
	\caption{The plots of the interpretations on $[0,1]$ of the terms in $T_{2}$.}
	\label{fig:piecewise}
\end{figure}
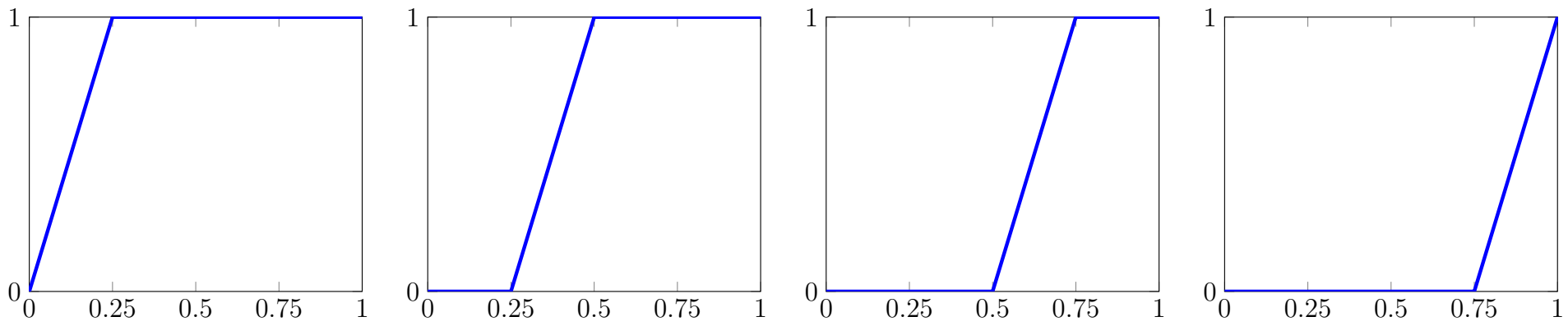

\begin{lemma} \label{l:piecewise-linear}
	For every $n \in \Nzero$, the interpretations on $[0,1]$ of the terms in $T_{n}$ are exactly the functions $f_{n,k} \colon [0,1]\to[0,1]$ of the form: 
    \[
    f_{n,k}(x) =	{
    	\begin{cases}
    			0				& \text{if } x \in \left[0, \frac{k}{2^n}\right];\\
    			2^n x - k			& \text{if } x \in \left[\frac{k}{2^n}, \frac{k+1}{2^n}\right];\\
    			1				& \text{if } x \in \left[\frac{k+1}{2^n}, 1\right];									\end{cases}
    	}
    \]
for  $k \in \{0,\dots,2^n-1\}$.
(See \cref{fig:piecewise} for plots of the functions \(f_{2,k}\) with $k \in \{0, 1, 2, 3\}$.)
\end{lemma}

\begin{proof}
Recall that, in the standard MV-algebra $[0,1]$,
\[
s(x)\df x\oplus x=\min\{1,2x\},\qquad p(x)\df x\odot x=\max\{0,2x-1\}.
\]
By the definition of $T_n$, for every $\tau\in T_n$, its interpretation $\tau^{[0,1]}\colon[0,1]\to[0,1]$ is a composition of $n$ maps, each equal to $s$ or $p$.
We prove the statement by induction on $n$.

\emph{Base case $n = 0$.}
The set $T_0$ is just the identity function, which is $f_{0,0}$.

\emph{Inductive step.}
Assume the claim holds for a fixed $n\in\Nzero$, and observe that every term in $T_{n+1}$ is either of the form $x\mapsto s(\sigma(x))$
or $x\mapsto p(\sigma(x))$ for some $\sigma\in T_n$. By the inductive hypothesis,
$\sigma^{[0,1]}=f_{n,k}$ for some $k\in\{0,\dots,2^n-1\}$. Hence it suffices to compute
$s\circ f_{n,k}$ and $p\circ f_{n,k}$.

\emph{Claim: ${s} \circ f_{n,k} = f_{n+1,2k}$.}
If $x\le \frac{k}{2^n}$ then $f_{n,k}(x)=0$ and thus $(s\circ f_{n,k})(x)=0$.
If $x\ge \frac{k+1}{2^n}$ then $f_{n,k}(x)=1$ and thus $(s\circ f_{n,k})(x)=1$.
On the middle interval $\left[\frac{k}{2^n},\frac{k+1}{2^n}\right]$ we have $f_{n,k}(x)=2^n x-k$.
Now $s(y)=\min\{1,2y\}$, and so $s(y)=2y \Leftrightarrow y\le \frac12$ and $s(y)=1 \Leftrightarrow y\ge \frac12$.
The equality $f_{n,k}(x)=\frac12$ holds exactly at $x=\frac{2k+1}{2^{n+1}}$.
Therefore:
\[
(s\circ f_{n,k})(x)=
\begin{cases}
0 & x\in\left[0,\frac{2k}{2^{n+1}}\right],\\[1mm]
2(2^n x-k)=2^{n+1}x-2k & x\in\left[\frac{2k}{2^{n+1}},\frac{2k+1}{2^{n+1}}\right],\\[1mm]
1 & x\in\left[\frac{2k+1}{2^{n+1}},1\right],
\end{cases}
\]
which is exactly $f_{n+1,2k}$.

One similarly obtains that $p\circ f_{n,k} = f_{n+1,2k+1}$.

We conclude that the interpretations of terms in $T_{n+1}$ are contained in
\[
\{f_{n+1,2k},f_{n+1,2k+1}\mid k=0,\dots,2^n-1\}=\{f_{n+1,j}\mid j=0,\dots,2^{n+1}-1\}.
\]
Conversely, every $f_{n+1,j}$ is of the form $s\circ f_{n, j/2}$ if $j$ is even, and
$p\circ f_{n,(j-1)/2}$ if $j$ is odd, hence it is the interpretation of a term in $T_{n+1}$.
This proves the equality of the two sets, completing the induction.
\end{proof}

\begin{lemma} \label{l:distance by tau}
Let $A$ be an MV-algebra, let $a,b \in A$, and let $n \in \Nzero$.
If for every $\tau \in T_{n}$ we have $\dMV( \tau(a),\tau(b))<1$, then
\[
\cha(a,b) \preceq \frac{1}{2^{n-1}}.
\]
\end{lemma}
\begin{proof}
By the definition of $\dMV$ (see \eqref{eq:norm-definition}), the claim is equivalent to:
\begin{equation} \label{eq:forallexists}
 \text{If } \forall \tau \in T_{n}\ \exists q_\tau \in [0,1) \cap \Q \text{ s.t.\ }\cha(\tau(a),\tau(b)) \preceq q_\tau, \text{ then }\cha(a,b) \preceq \frac{1}{2^{n-1}}.
\end{equation}

Note that, by the finiteness of $T_n$, the condition
\[
\forall \tau \in T_{n}\ \exists q_\tau \in [0,1) \cap \Q \text{ s.t.\ }\cha(\tau(a),\tau(b)) \preceq q_\tau
\]
appearing in \eqref{eq:forallexists}
is equivalent to
\[
\exists q \in [0,1) \cap \Q \ \forall \tau \in T_{n} \,\cha(\tau(a),\tau(b)) \preceq q.
\]
(Indeed, for the left-to-right direction, it suffices to take $q = \max_{\tau \in T_n} q_\tau$, and use  the property \eqref{eq:preceq-well-defined}.)

Therefore, \eqref{eq:forallexists} is equivalent to:
\begin{center}
$\forall q \in [0,1) \cap \Q$, if $\forall\tau \in T_{n}$, $\cha( \tau(a),\tau(b)) \preceq q$, then $\cha(a,b) \preceq \frac{1}{2^{n-1}}$.
\end{center}

This shows that the claim is equivalent to a set of quasi-equations. Thus, by \cref{t:dinola}, it is enough to check that they hold in $[0,1]$. 
Let $a,b \in[0,1]$ and assume without loss of generality that $a\leq b$. Suppose by contraposition that  $\cha(a, b)  \not\preceq  \frac{1}{2^{n-1}}$. In $[0,1]$ this is equivalent to $\lvert a -b\rvert>\frac{1}{2^{n-1}}$. Thus, there is $k \in \{0,\dots, 2^n-1\}$ such that $\left[\frac{k}{2^n}, \frac{k+1}{2^n}\right] \subseteq (a, b)$.
Let $f_{k}$ be the function defined in \Cref{l:piecewise-linear} for this particular $k$.  Notice that $f_{k}(a) = 0$ and $f_{k}(b) = 1$. By \cref{l:piecewise-linear}, there is $\tau \in T_{n}$ such that, for every $x \in[0,1]$, $f_k(x) = \tau^{[0,1]}(x)$. Hence,
\[
    \lvert \tau(a) - \tau(b) \rvert  =  \lvert f_{k}(a) - f_{k}(b) \rvert  =  \lvert 0 - 1 \rvert = 1,
\]
	which concludes the proof.
\end{proof}

\begin{lemma} \label{l:iterating tau}
Let $A$ be an MV-algebra endowed with an operation $\gamma$ of countably infinite arity satisfying \eqref{s:finite-axiom:item1}--\eqref{s:finite-axiom:item6}. 
If $\succ{a}$ is an MV-fast sequence in $A$, $n \in \Nzero$ and $\tau \in T_{n}$, then
\[\cha(\tau(\gamma(\asucc{a})), \tau(a_{n+1})) \preceq \frac{1}{2}.\]
\end{lemma}

\begin{proof}
Let $\tau \in T_{n}$. Iterating applications of \eqref{s:finite-axiom:item4} and \eqref{s:finite-axiom:item5}, one gets 
$\tau(\gamma(\asucc{a})) = \gamma(\tau(a_{n+1}),\tau(a_{n+2}),\dots)$; hence, 
\begin{align} \label{l:iterating tau:eq2}
\cha\mleft(\tau(\gamma(\asucc{a})), \tau(a_{n+1}) \mright) = 
\cha\mleft(\gamma(\tau(a_{n+1}),\tau(a_{n+2}),\dots),\tau(a_{n+1})\mright).
\end{align}
In addition, by \eqref{s:finite-axiom:item3},
\begin{align} \label{l:iterating tau:eq3}
\cha\mleft(\gamma(\tau(a_{n+1}),\tau(a_{n+2}),\dots),\tau(a_{n+1}) \mright) \preceq \frac{1}{2}.
\end{align}
Combining \eqref{l:iterating tau:eq2} and \eqref{l:iterating tau:eq3}, one gets the claim.
\end{proof}

\begin{theorem} \label{t:B1-B6-imply-CMV}
An MV-algebra endowed with an operation $\gamma$ of countably infinite arity is a $\CMV$-algebra if and only if it satisfies \eqref{s:finite-axiom:item1}--\eqref{s:finite-axiom:item6}.
\end{theorem}

\begin{proof}
The fact that every CMV-algebra satisfies \eqref{s:finite-axiom:item1}--\eqref{s:finite-axiom:item6} is the content of \cref{p:soundness finite}. 
For the other direction,
first notice that \eqref{d:CMV:A1} and \eqref{d:CMV:A2} coincide with \eqref{s:finite-axiom:item1} and \eqref{s:finite-axiom:item2}. Thus, $A$ is a semisimple by \Cref{p:CMV arch} and hence we are left to prove \eqref{d:CMV:A3b}.
Consider an arbitrary sequence $\succ{a}$ in $A$. By \eqref{s:finite-axiom:item6},
\begin{align} \label{eq:proof:C3.1}
\cha(\gamma(\asucc{a}),\rho_{n}(\init{a})) = \cha(\gamma(\asucc{\rho(a)}),\rho_{n}(\init{a})).
\end{align}
By \Cref{l:in V dist rho_n rho_m}, $\succ{\rho(a)}$ is MV-fast; hence, we can invoke \Cref{l:iterating tau} to obtain, for every $n \in \Npos$ and $\tau \in T_{n-1}$, 
\[\cha\mleft(\tau(\gamma(\asucc{\rho(a)})), \tau(\rho_{n}(\init{a})) \mright) \preceq \frac{1}{2}.\] 
It follows from \Cref{l:leqisleq} that:
\[\dMV{\mleft(\tau(\gamma(\asucc{\rho(a)})), \tau(\rho_{n}(\init{a})) \mright)}=\normMV{\cha\mleft(\tau(\gamma(\asucc{\rho(a)})), \tau(\rho_{n}(\init{a})) \mright)}\leq \frac{1}{2}.\] 
Thus, by \Cref{l:distance by tau},
\begin{align} \label{eq:proof:C3.2}
\cha(\gamma(\asucc{\rho(a)}),\rho_{n}(\init{a})) \preceq \frac{1}{2^{n-2}}.
\end{align}
Combining \eqref{eq:proof:C3.1} and \eqref{eq:proof:C3.2} we see that $\cha(\gamma(\asucc{a}), \rho_{n}(\init{a})) \preceq \frac{1}{2^{n-2}}$.  Therefore, $\rho_{n}(\init{a})$ converges to $\gamma(\asucc{a})$.  By \Cref{l:Cauchy}, $\dMV(\gamma(\asucc{a}),\rho_{n}(\init{a})) \leq \frac{1}{2^n}$.  Therefore, axiom \eqref{d:CMV:A3b} holds. So, by \Cref{r:eqA3} $A$ is a $\CMV$-algebra.
\end{proof}

\section{Free metrically complete MV-algebras}

Having proved that $\cMV$ is equivalent to a variety, it makes sense to call \emph{free} a metrically complete MV-algebra that, through the equivalence, corresponds to a free $\CMV$-algebra.

In this section, we provide a concrete description of the free $\CMV$-algebra on a set $I$ as an algebra of $[0,1]$-valued functions on the power $[0,1]^I$. As a consequence, we obtain that $\CMV$ is the category of algebras of the (infinitary) Linton--Lawvere theory whose objects are powers of $[0,1]$ and whose morphisms are the continuous maps that \emph{decrease denominators} (see definitions below).

We begin with a purely algebraic result of independent interest.
We recall the standard operators $\mathbb{I}$ (closure under isomorphisms), $\mathbb{S}$ (closure under subalgebras) and $\mathbb{P}$ (closure under products).
Moreover, we denote simply with $[0,1]$ the canonical $\CMV$-algebra whose underlying set is the unit interval $[0,1]$. 
\begin{theorem} \label{t:W-equal-ISP}$
\CMV=\ISP([0,1])$.
\end{theorem}
\begin{proof}
The inclusion $\supseteq$ is clear because $\CMV$ is a variety containing $[0,1]$. For the inclusion $\subseteq$, let $A \in \CMV$. Then $A$ is a semisimple MV-algebra, and thus the morphisms of MV-algebras towards $[0,1]$ separate the points of $A$. Since the forgetful functor $\U \colon \CMV\to \MV$ is full, every morphism of MV-algebras from $A$ to $[0,1]$ is also a $\CMV$-morphism. Hence, there are enough $\CMV$-morphisms from $A$ to $[0,1]$ to separate the points of $A$, and so $A$ is isomorphic to a subalgebra of a power of $[0,1]$.
\end{proof}

By \cref{t:W-equal-ISP}, the operator sending an equivalence class of $\CMV$-terms (where two terms are equivalent if their interpretation in each algebra of the variety coincides)
to the interpretation in $[0,1]$ is an injective operator.
We next characterise its image, as the class of so-called \emph{denominator-decreasing continuous maps}, also called \emph{a-maps}.

Recall from \cite[Def.~2.2]{AbbadiniMarraEtAl2025} that an \emph{a-space} is a topological space $X$ equipped with a map $\zeta \colon X \to \Nzero$ (called the \emph{denominator map})  that is continuous w.r.t.\ the upper topology of the divisibility lattice of $\Nzero$, i.e., for every $n \in \Nzero$, the preimage under $\zeta$ of the set of divisors of $n$ in $\Nzero$ is closed.  A map $f \colon X \to Y$ between a-spaces is called \emph{denominator-decreasing}, if, for every $x \in X$, $\zeta_Y(f(x))$ divides $\zeta_X(x)$. A continuous denominator-decreasing map is called an \emph{a-map}. 
Examples of a-spaces are $[0,1]$ with denominator map given by:
\[
\den(x)\coloneqq\begin{cases}
q & \text{if } x=\frac{p}{q} \text{ with $p$ and $q$ coprime}\\
0 & \text{if } x\notin \Q,
\end{cases}
\]
and $[0,1]^I$ with denominator map given by $\den_I(x)\coloneqq\lcm\{\den x_i\mid i\in I\}$, where we recall that ``$\lcm$'' means ``least upper bound in the complete lattice $\Nzero$ with respect to the divisibility order'' (which has $1$ as minimum and $0$ as maximum).

\begin{lemma} \label{l:gamma is a-map}
The function $\gamma \colon [0,1]^{\Npos}\to [0,1]$ is an a-map. 
\end{lemma}
\begin{proof}
For any $k \in \Npos$, let us define a function $\tilde{\rho}_k \colon [0,1]^{\Npos}\to [0,1]$ by setting
\begin{align*}
\tilde{\rho}_k\left(\succ{x}\right) \df \rho_k(x_1,\dots,x_{k}).
\end{align*}
For every $i \in \Npos$, the projection $\pi_i \colon [0,1]^{\Npos}\to [0,1]$ onto the $i$-th coordinate is continuous, and, for every $n \in \Npos$, $\rho_n \colon [0,1]^n\to [0,1]$ is continuous. Therefore, for every $k \in \Npos$, $\tilde{\rho}_k \colon [0,1]^\Npos\to [0,1]$ is continuous. The sequence  $\succ{\tilde{\rho}}$ of continuous functions converges to $\gamma$ uniformly, because
\[
\sup_x\left|\gamma(x)-\widetilde\rho_k(x)\right|\leq 2^{-k},
\]
which follows from \eqref{d:CMV:A3b}.
Therefore, the uniform limit 
$\gamma$ is continuous. Moreover, for each $n \in \Npos$, $\rho_n$ decreases denominators, and thus $\tilde{\rho}_n$, as well. Since the sequence $\succ{\tilde{\rho}}$ of denominator-decreasing functions converges pointwise to $\gamma$, $\gamma$ decreases denominators; indeed, for every $x \in [0,1]^\Npos$, if $\den_\Npos(x)=d \neq 0$ then the possible values in $[0,1]$ whose denominator divides $d$ form a finite closed set and so $\lim_{k \to \infty}\widetilde\rho_k(x)$ stays in that set.
\end{proof}

The following theorem is an easy adaptation of \cite[Thm.~7.7]{AbbadiniMarraEtAl2025}, which was in the context of unital Abelian $\ell$-groups. 

\begin{theorem}[Stone-Weierstrass for MV-algebras]\label{t:StWe}
    Let $X$ be a compact a-space and $A$ an MV-subalgebra of the MV-algebra of a-maps from $X$ to $[0,1]$ (with pointwise defined operations) such that
    \begin{enumerate}
        \item \label{i:sep} for every $x \neq y \in X$, there is $f \in A$ such that $f(x) \neq f(y)$;
        \item \label{i:pres-den} for every $x \in X$, $\zeta(x) = \lcm\{\den(f(x))\mid f \in A\}$.
    \end{enumerate} 
    Then $A$ is dense in the set of a-maps from $X$ to $[0,1]$ with respect to the uniform metric.
\end{theorem}

The primitive function symbols of $\CMV$ (i.e., $\gamma$, $\oplus$, $\lnot$, $0$) have a natural interpretation in $[0,1]$; each of them is an a-map from a power of $[0,1]$ to $[0,1]$ itself. In fact, a-maps characterise terms, in the following sense:

\begin{theorem} \label{t:operations-coincide-with-continuous-den-decreasing}
For every set $I$, the set of a-maps from $[0,1]^I$ to $[0,1]$ coincides with the set of interpretations in $[0,1]$ of the $\CMV$-terms of arity $I$.
\end{theorem}
\begin{proof}
Let $T$ be the set of interpretations in $[0,1]$ of the $\CMV$-terms of arity $I$ and let $A$ be the set of a-maps from $[0,1]^I$ to $[0,1]$.  It easily follows from \Cref{l:gamma is a-map} that $T$ is a sub-MV-algebra of $A$. To prove the equality, notice that, for each $i \in I$, $T$ contains the projection $\pi_i \colon [0,1]^I\rightarrow [0,1]$. Thus, (i) the maps in $T$ separate the points, i.e.\ for every $x\neq y\in[0,1]^I$ there is $g\in T$ such that $g(x)\neq g(y)$; and (ii) for every $x\in [0,1]^I$, $\den_I{x}=\lcm\{{\den{\pi_i(x)}}\}_{i\in I}$.  Thus, \Cref{t:StWe} applies, and we obtain that $T$ is dense in $A$. Furthermore, $T$ is a $\CMV$-algebra, and therefore $T$ is Cauchy complete; thus $T=A$.
\end{proof}

Roughly speaking, by \cref{t:W-equal-ISP}, sending an equivalence class of $\CMV$-terms (where two terms are equivalent if their interpretation in each algebra of the variety coincides)
to the interpretation in $[0,1]$ is an injective operator, and, by
\cref{t:operations-coincide-with-continuous-den-decreasing}, it is surjective, too.
So, the equivalence classes of $\CMV$-terms of arity $I$ are in bijective correspondence with the a-maps from $[0,1]^I$ to $[0,1]$.

\begin{theorem}\label{t:free-CMV-algebra}
Let $I$ be a set. The $\CMV$-algebra 
\[
C\df \{f \colon[0,1]^I\to [0,1] \mid f \text{ is an a-map}\}
\]
is freely generated by the projections $\left(\pi_i \colon [0,1]^I\to [0,1]\right)_{i \in I}$.
\end{theorem}
\begin{proof}
By \Cref{t:W-equal-ISP}, the free algebra generated by the set $I$ is the set of functions $[0,1]^I\to [0,1]$ that are the interpretation in $[0,1]$ of an operation of arity $I$. By \Cref{t:operations-coincide-with-continuous-den-decreasing}, this set is precisely $C$.
\end{proof}

For the reader familiar with Lawvere-Linton's theories (\cite{Lawvere,Linton}), we mention that the results in this section show that $\CMV$ is the category of algebras of the varietal theory of a-maps between powers of $[0,1]$.

\section{Metrically complete MV-algebras are not equivalent to a finitary variety}

In this section, we prove that the equivalent categories $\cG$, $\cMV$ and $\CMV$ are not equivalent to any \emph{finitary} variety. In fact, we prove more: they are not equivalent to any elementary class of finitary structures, with morphisms all the homomorphisms.

Recall that, in elementary classes, directed colimits are preserved by the forgetful functor to the category $\Set$ \cite{Richter1971}.

Let $\KH$ be the category of compact Hausdorff spaces and continuous functions.
We will use the following result.

\begin{theorem} \label{t:no-faithful}
    No faithful functor from $\KH^\op$ to $\Set$ preserves directed colimits.
\end{theorem}

\begin{proof}
    At the last line of the proof of \cite[Thm.~22]{LiebermanRosickyEtAl2023}, the authors conclude that no faithful functor from the category of commutative unital $C^*$-algebras to $\Set$ preserves directed colimits.
    The result in the statement follows since, by Gelfand duality, the category of commutative unital $C^*$-algebras is equivalent to $\KH^\op$.
\end{proof}
For a short and direct proof of \cref{t:no-faithful}, that does not pass via Gelfand duality, we refer to \cite{Abbadini2026}.

As observed in \cite[Cor.~24]{LiebermanRosickyEtAl2023}, it follows that $\KH^{\rm op}$ is not equivalent to any elementary (i.e., first-order definable) class of structures with homomorphisms as morphisms. We now show that \cref{t:no-faithful} implies that $\CMV$ is not equivalent to any elementary class of structures.
Recall from above that an a-space is a topological space $X$ equipped with a map $\zeta \colon X \to \Nzero$ such that, for every $n \in \Nzero$, the preimage under $\zeta$ of the set of divisors of $n$ in $\Nzero$ is closed.
Recall from \cite[Def.~4.1]{AbbadiniMarraEtAl2025}
that a \emph{normal a-space} is an a-space such that:
\begin{enumerate}
\item (Compactness) $X$ is compact;
\item (Separation) for any two distinct points $x,y \in X$, there are disjoint open neighbourhoods
$U$ and $V$ of $x$ and $y$, respectively, such that $\zeta(z)=0$ for every
$z \in X \setminus (U \cup V)$.
\end{enumerate}
We denote by $\KHa$ the category of normal a-spaces and a-maps.
The main result in \cite{AbbadiniMarraEtAl2025} is the following.
\begin{theorem}[{\cite[Thm.~7.10]{AbbadiniMarraEtAl2025}}] \label{t:metr-complete-duality}
    The category $\cG$ of metrically complete $\ell$-groups is dually equivalent to the category $\KHa$ of normal a-spaces.
\end{theorem}

    There is an adjunction
    \[
    \begin{tikzcd}
    \KHa \arrow[yshift=5pt]{rr} & {\adj} & \KH \arrow[yshift=-5pt]{ll}{\Delta},
    \end{tikzcd}
    \]
    where the functor $\KHa \to \KH$ is the forgetful functor, and
    \[
    \Delta \colon \KH \to \KHa
    \]
    sends a compact Hausdorff space $X$ to the same space equipped with the
    denominator function constantly equal to $0$, and sends a continuous map to itself.
    \begin{lemma}\label{l:Delta-preserves-limits}
    The functor $\Delta$ preserves limits over non-empty diagrams.
    \end{lemma}
    \begin{proof}
    If $C$ is the limit of a diagram $\mathcal{D} \colon \cat{I} \to \KH$ in $\KH$ with $\mathsf{I}$ non-empty, then $\Delta C$ is a limit of $\Delta \mathcal{D}$ because the vertex of every cone in $\KHa$ over $\Delta \mathcal{D}$ must have a denominator map that evaluates all the points to $0$. This is because there is at least an object in the image of $\Delta\mathcal{D}$ and all its points have denominator $0$. Thus, the domain of every a-map into it must also have all points of denominator $0$. It follows that the vertex $V$ of any cone over $\Delta \mathcal{D}$ has all points of denominator $0$.
    Therefore, the canonical map from $V$ to the vertex of $C$ decreases denominators.    
    \end{proof}

\begin{theorem} \label{t:faithful}
    No faithful functor from $\KHa^\op$ to $\Set$ preserves directed colimits.
\end{theorem}

\begin{proof}
    By \cref{l:Delta-preserves-limits},
    \[
    \Delta^\op \colon \KH^\op \to \KHa^\op
    \]
    preserves directed colimits (necessarily non-empty); moreover, it is clearly faithful.
    Therefore, if there were a faithful functor $F\colon \KHa^\op \to \Set$ preserving directed colimits, then the composite 
    \[
    \KH^\op \xlongrightarrow{\Delta^\op} \KHa^\op \xlongrightarrow{F} \Set
    \]
    would be a faithful functor that preserves directed colimits, contradicting \cref{t:no-faithful}.
\end{proof}

\begin{corollary}\label{c:no-preserves}
    No faithful functor from the category $\cG$ of metrically complete $\ell$-groups  to $\Set$ preserves directed colimits.
\end{corollary}
\begin{proof}
    Since $\cG$ is dually equivalent to $\KHa$ by \cref{t:metr-complete-duality}, the result follows from \cref{t:faithful}.
\end{proof}

\begin{corollary} \label{c:not-elementary}
  The categories $\cMV$ and $\cG$ of metrically complete MV-algebras and metrically complete $\ell$-groups are not equivalent to any elementary class of finitary structures (nor, in particular, to any variety of finitary algebras).
\end{corollary}

\begin{proof}
    This follows from \cref{c:no-preserves} and the fact that, in elementary classes, directed colimits are preserved by the forgetful functor to the category $\Set$ of sets and functions~\cite{Richter1971}.
\end{proof}

Finally, we use \cref{c:not-elementary} to prove a further negative result. By an \emph{$\ISP$-class} of algebras we mean a class of algebras closed under isomorphisms, subalgebras and arbitrary products.
When we consider it as a category, the morphisms are the homomorphisms. In such a category,  closure under homomorphic images can be expressed by a purely categorical property, namely, \emph{Barr-exactness}\footnote{We recall that a regular category $\cat{C}$ (such as any class of algebras closed under isomorphisms, subalgebras and finite products) is \emph{Barr-exact} if every equivalence relation is the kernel pair of some morphism. Roughly speaking, this says that the quotient under a congruence on an algebra in the class still lies in the class.}.

Now, it is clear that Barr-exactness is invariant under categorical equivalence. This leads to the following result.
\begin{corollary} \label{c:SP-class}
  The categories $\cMV$ and $\cG$ of metrically complete MV-algebras and metrically complete $\ell$-groups  are not equivalent to any $\ISP$-class of finitary algebras.
\end{corollary}
\begin{proof}
    Since $\cMV$ is equivalent to a variety of infinitary algebras (\cref{t:cG-equiv-CMV}), it is Barr-exact.
  If an $\ISP$-class of finitary algebras is Barr-exact as a category, then it is closed under homomorphic images, and so it is a variety of finitary algebras (see \cite{Banaschewski83}).
    By \cref{c:not-elementary}, $\cMV$ is not equivalent to any variety of finitary algebras.
  Therefore, $\cMV$ is not equivalent to any $\ISP$-class of finitary algebras.
\end{proof}
The statement of \cref{c:SP-class} involves finitary algebras only. However, we emphasise in closing that our proof does rely on our main result \cref{t:cG-equiv-CMV} that the category $\cMV$ of metrically complete MV-algebras is equivalent to some variety of infinitary algebras.

\bibliography{AMS_Biblio2}
\bibliographystyle{plain}

\end{document}